\documentclass{article}

\usepackage{microtype}
\usepackage{graphicx}
\usepackage{booktabs} 

\usepackage{hyperref}

\usepackage[accepted]{icml2020}

\usepackage{amsmath}
\usepackage{amssymb}
\usepackage{amsfonts}
\usepackage{amsthm}
\usepackage{mathtools}

\usepackage[utf8]{inputenc}
\usepackage[T1]{fontenc}

\DeclarePairedDelimiter{\ceil}{\lceil}{\rceil}

\usepackage[sans]{dsfont}

\usepackage{acronym}
\usepackage{latexsym}
\usepackage{xspace}

\newcommand{\dd}{\:d}

\newcommand{\eps}{\varepsilon}

\newcommand{\pd}{\partial}

\newcommand{\dif}{\dd}

\DeclareMathOperator*{\argmin}{argmin}

\DeclareMathOperator{\diam}{diam}

\DeclareMathOperator{\dom}{dom}

\DeclareMathOperator{\tr}{tr}

\newcommand{\ca}{\mathtt{a}}
\newcommand{\cb}{\mathtt{b}}
\newcommand{\ce}{\mathtt{e}}
\newcommand{\ct}{\mathtt{t}}

\newcommand{\bE}{\mathbf{E}}

\renewcommand{\emptyset}{\varnothing}
\newcommand{\eqdef}{\triangleq}

\newcommand{\scrA}{\mathcal{A}}

\newcommand{\scrF}{\mathcal{F}}

\newcommand{\scrK}{\mathcal{K}}
\newcommand{\scrL}{\mathcal{L}}

\newcommand{\scrO}{\mathcal{O}}
\newcommand{\scrP}{\mathcal{P}}

\newcommand{\scrS}{\mathcal{S}}

\newcommand{\scrW}{\mathcal{W}}
\newcommand{\scrX}{\mathcal{X}}

\DeclareMathOperator{\mat}{mat}

\newcommand{\0}{\mathbf{0}}

\newcommand{\Rn}{\R^n}

\newcommand{\R}{\mathbb{R}}

\newcommand{\N}{\mathbb{N}}

\DeclareMathOperator{\ball}{\mathbb{B}}

\newcommand{\bC}{{\mathbf{C}}}

\newcommand{\metric}{\mathsf{d}}

\DeclareMathOperator{\gap}{\mathsf{Gap}}

\theoremstyle{plain}
\newtheorem{theorem}{Theorem}

\newtheorem*{corollary*}{Corollary}
\newtheorem{lemma}[theorem]{Lemma}
\newtheorem{proposition}[theorem]{Proposition}

\theoremstyle{definition}
\newtheorem{definition}[theorem]{Definition}
\newtheorem*{definition*}{Definition}

\theoremstyle{remark}

\newtheorem*{remark*}{Remark}
\newtheorem*{notation*}{Notational remark}
\newtheorem{example}{Example}

\numberwithin{theorem}{section}
\numberwithin{remark}{section}
\numberwithin{example}{section}

\DeclarePairedDelimiter{\abs}{\lvert}{\rvert}

\DeclarePairedDelimiter{\inner}{\langle}{\rangle}
\DeclarePairedDelimiter{\norm}{\lVert}{\rVert}

\newacro{HBA}[HBA]{Hessian-barrier algorithm}
\newacro{AHBA}[AHBA]{Adaptive-Hessian-barrier algorithm}
\newacroplural{NE}[NE]{Nash equilibria}
\newacro{VI}{variational inequality}
\newacroplural{VI}[VIs]{variational inequalities}
\newacro{iid}[i.i.d.]{independent and identically distributed}
\newacro{FOM}{First-order method}
\newacro{FW}{Frank-Wolfe}
\newacro{CG}{Conditional Gradient}
\newacro{LMO}{Linear minimization oracle}
\newacro{LLOO}{Local Linear minimization oracle}
\newacro{SC}{self-concordant}

\newcommand{\subf}[2]{%
  {\small\begin{tabular}[t]{@{}c@{}}
  #1\\#2
  \end{tabular}}%
}

\def\pd#1{{\color{blue}#1}}

\icmltitlerunning{Self-concordant analysis of Frank-Wolfe algorithms}
\begin{document}

\twocolumn[
\icmltitle{Self-Concordant Analysis of Frank-Wolfe Algorithms}



\icmlsetsymbol{equal}{*}

\begin{icmlauthorlist}
\icmlauthor{Pavel Dvurechensky}{WIAS,HSE}
\icmlauthor{Petr Ostroukhov}{MI}
\icmlauthor{Kamil Safin}{MI}
\icmlauthor{Shimrit Shtern}{TEC}
\icmlauthor{Mathias Staudigl}{DKE}
\end{icmlauthorlist}

\icmlaffiliation{DKE}{Department of Data Science and Knowledge Engineering, Maastricht University, The Netherlands,}
\icmlaffiliation{WIAS}{Weierstrass Institute for Applied Analysis and Stochastics, Berlin, Germany.}
\icmlaffiliation{TEC}{Grand Fellow, Technion - Israel Institute of Technology, Haifa, Israel.}
\icmlaffiliation{MI}{Moscow Institute of Physics and Technology, Dolgoprudny, Russia}
\icmlaffiliation{HSE}{HSE University, Moscow, Russia}
\icmlcorrespondingauthor{Mathias Staudigl}{mathias.staudigl@gmail.com}

\icmlkeywords{Self-concordant function, Frank-Wolfe method, Conditional Gradient method, linear minimization oracle, linear convergence rates}

\vskip 0.3in
]



\printAffiliationsAndNotice{\icmlEqualContribution} 

\begin{abstract}
Projection-free optimization via different variants of the \ac{FW}, a.k.a. Conditional Gradient method has become one of the cornerstones in optimization for machine learning since in many cases the linear minimization oracle is much cheaper to implement than projections and some sparsity needs to be preserved. In a number of applications, e.g. Poisson inverse problems or quantum state tomography, the loss is given by a \ac{SC} function having unbounded curvature, implying absence of theoretical guarantees for the existing \ac{FW} methods. We use the theory of \ac{SC} functions to provide a new adaptive step size for \ac{FW} methods and prove global convergence rate $O(1/k)$ after $k$ iterations.
If the problem admits a stronger local linear minimization oracle, we construct a novel \ac{FW} method with linear convergence rate for \ac{SC} functions.
\end{abstract}

\section{Introduction}
\label{sec:Introduction}
Statistical analysis using (generalized) self-concordant (SC) functions as a loss function is gaining increasing attention in the machine learning community \cite{Bac10,Owe13,OstBac18,MarBacOstRu19}.  This class of loss functions allows to obtain faster statistical rates akin to least-squares. At the same time, the minimization of empirical risk in this setting is a challenging optimization problem in high dimensions. The reason is that for large-scale problems usually a \ac{FOM} is the method of choice. Yet, self-concordant functions are usually not strongly convex and do not have Lipschitz continuous gradients. Hence, unless we can verify a relative smoothness condition \cite{BauBolTeb17,LuFreNes18},  the cornerstone assumptions for \ac{FOM} do not apply when minimizing \ac{SC} functions. Further, many machine learning problems involve some sparsity constraint usually admitting a \ac{LMO} \cite{Jag13}, meaning that for a reasonable computational price one can minimize a linear function over the feasible set; Examples include $\ell_{1}$ regularization, Spectrahedron constraints \cite{Haz08}, and ordered weighted $\ell_1$ regularization \cite{FigZen14}. In such settings an attractive method of choice is the Frank-Wolfe (\ac{FW}) method \cite{FraWol56}, also known as the \ac{CG}. The modern convergence analysis of \ac{FW}, starting with \cite{Jag13}, relies on bounded curvature of the objective function $f$. A sufficient condition for this is that the objective function has a Lipschitz continuous gradient over the domain of interest. In general, \ac{SC} functions have unbounded curvature and are not strongly convex. In fact, to the best of our knowledge, there is no existing convergence guarantee for \ac{FW} when minimizing a general \ac{SC} function. Given the plethora of examples of \ac{SC} loss functions, and the importance of \ac{FW} in machine learning, we focus on the question of convergence and complexity of \ac{FW} for minimizing a SC function $f$ over a set $\scrX$ with \ac{LMO}. 

\textbf{Our Contributions. } In this paper we demonstrate that \ac{FW} indeed works when minimizing a \ac{SC} function over a compact convex domain. However, since the standard arguments for \ac{FW} are not applicable, we have to develop some new ideas. Section \ref{sec:sublinear} constructs new adaptive variants of the \ac{FW} method, whose main innovation is the construction of new step-size policies ensuring global convergence and $\scrO(1/k)$ convergence rates. One of the keys to our results is to develop step-sizes which ensure monotonicity of the method, leading to Lipschitz smoothness and strong convexity of the objective on a level set. While the Lipschitz smoothness on a level set is used in the analysis, the algorithms do not require knowledge of the Lipschitz constant but rather rely only on basic properties of \ac{SC} functions.
Section \ref{sec:fast} gives a \emph{linear} convergence result for \ac{FW}, under the assumption that we have access to the strong convexity parameter on the level set and a \ac{LLOO}, first constructed in \citet{GarHaz16}. 

\textbf{Previous work. } Despite its great scalability properties, \ac{FW} is plagued by slow convergence rates. \citet{FraWol56} showed an $\scrO(1/k)$ convergence rate of the function values to the optimal value for quadratic programming problems. Later on \citet{LevPol66} generalized this result to general functions and arbitrary compact convex sets, and \citet{CanCul68} proved that this rate is actually tight (see also \citet{Lan13}). However, a standard result in convex optimization states that projected gradient descent converges linearly to the minimizer of a strongly convex and Lipschitz (i.e. well-conditioned) smooth function. Subsequently, an important line of research emerged with the aim to accelerate \ac{FW}. \citet{GueMar86} obtained linear convergence rates in well conditioned problems over polytopes under the a-priori assumption that the solution lies in the relative interior of the feasible set, and the rate of convergence explicitly depends on the distance the solution is away from the boundary (see also \citet{EpeFreu00,BecTeb04}). If no a-priori information on the location of the solution is available, there are essentially two known twists of the vanilla \ac{FW} to boost the convergence rates. One twist is to modify the search directions via \emph{corrective} or \emph{away} search directions \cite{JagLac15,FreGriMaz17,BecSht17,GutPen20,PenRod18}. These methods require however more fine grained oracles which are not in general available for \ac{SC} functions\footnote{The away step version of \ac{FW} needs a vertex oracle. Vertices need not be in the domain of a \ac{SC} function (i.e. the origin in Example \ref{example}).}. The alternative twist is to restrict the \ac{LMO} in order to gain a more powerful local approximating model of the objective function. This last strategy has been used in several recent papers \cite{HarJudNem15,GarHaz16}. A unified analysis for different settings was recently proposed in \cite{PedNegAskJag20}, which inspired one of the variants of our methods.

None of these references explicitly studied the important case of \ac{SC} minimization. 
\citet{Nes18CG} requires the gradient of the objective function to be H\"older continuous. Implicitly it is assumed that the feasible set $\scrX$ is contained in the domain of the objective function $f$, an assumption we do not make, and is also not satisfied in important applications (e.g. $0\in\scrX$ in the Poisson inverse problem below, but $0\notin\dom f$). Specialized to solving a quadratic Poisson inverse problem in phase retrieval, \citet{Odor16} provided a globally convergent \ac{FW} method. They analyzed this problem via techniques developed in \cite{Nes18CG} and with step-size specific to the application. We instead give a unified general analysis valid for \emph{all} \ac{SC} functions. Given the importance of \ac{SC} minimization in machine learning, a large volume of work around proximal Newton and proximal gradient methods have been developed in a series of papers by \cite{CevKyrTra13,CevKyrTra14,CevKyrTra15}. This paper complements this line of research by providing a \ac{FW} analysis of \ac{SC} optimization problems. Recently, \cite{LiuCevTra20} coupled a \ac{LMO} with a second-order model of the \ac{SC} objective function and obtained similar results to ours. Similar to \cite{LiuCevTra20}, the recent paper \cite{carderera2020second-order} studied a Newton-\ac{FW} approach, but for well-conditioned functions. 

\section{Preliminaries}
\label{sec:prelims}
We consider the minimization problem
\begin{equation}\label{eq:P}
\min_{x\in\scrX}f(x),
\end{equation}
where $\scrX$ is a compact convex set in $\Rn$ with nonempty interior. We assume that $\dom f\cap\scrX\neq\emptyset$ and the set $\scrX$ is represented by an \ac{LMO} returning at a point $x$ the target vector 
\begin{equation}\label{eq:s}
s(x)\eqdef\arg\min_{s\in\scrX}\inner{\nabla f(x),s}.
\end{equation} 
In case of multiple solutions of \eqref{eq:s} we assume that ties are broken by some arbitrary mechanism. The function $f$ is assumed to be self-concordant (\ac{SC}) \cite{NesIPM94} in the following sense: 
Let $f\in\bC^{3}(\dom f)$ be a closed convex function with open domain $\dom f\eqdef\{x\in\Rn\vert f(x)<\infty\}$.  For a fixed $x\in\dom f$ and direction $u\in\Rn$, define $\phi(x;t)\eqdef f(x+tu)$ for $t\in\dom\phi(x;\cdot)$.
\begin{definition}\label{def:SC}
A proper closed convex function $f:\Rn\to(-\infty,\infty]$ with open domain $\dom f$ is self-concordant (SC) with parameter $M>0$ (i.e. $f\in\scrF_{M}$) iff 
\begin{equation}
\abs{\phi'''(x,0)}\leq M\phi''(x,0)^{3/2}.
\end{equation}
\end{definition}
\ac{SC} functions were originally introduced in the context of interior-point methods \cite{NesIPM94}, but recently have received significant attention in machine learning as many problems in supervised learning and empirical risk minimization involve loss functions which are \ac{SC}. We give a (non-exhaustive) list of examples next. Note that, unlike classical settings, $f$ in \eqref{eq:P} is not assumed to be a self-concordant barrier for $\scrX$.

\subsection{Examples}
\label{sec:examples}

\textbf{Poisson Inverse Problem. }
Consider the low-light imaging problem in signal processing, where the imaging data is collected by counting photons hitting a detector over time. In this setting, we wish to accurately reconstruct an image in low-light, which leads to noisy measurements due to low photon count. Assume that the data-generating process of the observations follows as Poisson distribution 
\[
p(y\vert Wx^{\ast})=\prod_{i=1}^{m}\frac{(w_{i}^{\top}x^{\ast})^{y_{i}}}{y_{i}!}\exp(-w_{i}^{\top}x^{\ast}),
\]
where $x^{\ast}\in\Rn$ is the true image, $W$ is a linear operator that projects the scene onto observations, $w_{i}$ is the $i$-th row of $W$ and $y\in\N^{m}$ is the vector of observed photon counts. The maximum likelihood formulation of the Poisson inverse problem \cite{HarMarWil11} under sparsity constraints leads to the optimization problem
\begin{equation}\label{eq:Poisson}
\min_{x\in\Rn_{+};\norm{x}_{1}\leq M} \{f(x)=\sum_{i=1}^{m}w_{i}^{\top}x-\sum_{i=1}^{m}y_{i}\ln(w_{i}^{\top}x)\}.
\end{equation}
Setting $\varphi_{i}(t)\eqdef t-y_{i}\ln(t)$ for $t>0$, we observe that $f(x)= \sum_{i=1}^{m}\varphi_{i}(w_{i}^{\top}\theta)$. Since each individual function $\varphi_{i}$ is \ac{SC} with parameter $M_{\varphi_{i}}\eqdef \frac{2}{\sqrt{y_{i}}}$, general rules of \ac{SC} calculus shows that $f$ is \ac{SC} with domain $\bigcap_{i=1}^{m}\{x\in\Rn\vert w_{i}^{\top}x>0\},$ and parameter $M\eqdef\max_{1\leq i\leq m}\frac{2}{\sqrt{y_{i}}}$ (\citet{Nes18}, Thm. 5.1.1).

\textbf{Learning Gaussian Markov random fields. }
We consider learning a Gaussian graphical random field of $p$ nodes/variables from a data set $\{\phi_{1},\ldots,\phi_{N}\}$ \cite{SpeKli86,CevKyrTra13}. Each random vector $\phi_{j}$ is an i.i.d realization from a $p$-dimensional Gaussian distribution with mean $\mu$ and covariance matrix $\Sigma$. Let $\Theta=\Sigma^{-1}$ be the precision matrix. To satisfy conditional dependencies between the random variables, $\Theta$ must have zero in $\Theta_{ij}$ if $i$ and $j$ are not connected in the underlying dependency graph. To learn the graphical model subject to sparsity constraints, we minimize the loss function 
$$
f(x)\eqdef -\ln\det(\mat(x))+\tr(\hat{\Sigma}\mat(x))
$$
over the $\ell_{1}$-ball $\scrX\eqdef \{x\in\Rn\vert \mat(x)\succ 0,\norm{x}_{1}\leq R\}$, where $n=p^{2}$ and $\mat(x)\in\R^{p\times p}_{\text{sym}}$ is the $p\times p$ symmetric matrix constructed from the $n=p^{2}$-dimensional vector $x$. $f$ is \ac{SC} with $M=2$.

\textbf{Logistic Regression. }
Consider the regularized logistic regression problem 
\begin{equation}\label{eq:Logistic}
f(x)=\frac{1}{N}\sum_{i=1}^{N}\ell(y_{i}(\inner{\phi_{i},x}+\mu))+\frac{\gamma}{2}\norm{x}^{2}_{2}
\end{equation}
where $\ell(t)\eqdef \log(1+e^{-t})$ is the logistic loss, $\mu$ is a given intercept, $y_{i}\in\{-1,1\}$ and $\phi_{i}\in\Rn$ are given as input data for $i=1,2,\ldots,N$, and $\gamma>0$ is a given regularization parameter. According to Prop. 5 in \citet{SunTran18}, this functions is \ac{SC} with parameter $M=\frac{1}{\sqrt{\gamma}}\max\{\norm{\phi_{i}}_{2}\vert 1\leq i\leq n\}$. To promote sparsity, we minimize \eqref{eq:Logistic} over the $\ell_{1}$-ball $\scrX=\{x\in\Rn\vert \norm{x}_{1}\leq R\}$.  The resulting \ac{SC} optimization problem \eqref{eq:P} is the constrained formulation of the elastic-net  regularization. 

\textbf{Portfolio Optimization. } In this problem there are $n$ assets with returns $r_{t}\in\R^{n}_{+}$ in period $t$ of the investment horizon. The goal is to minimize the utility function $f(x)$ of the investor by choosing the weights of the assets in the portfolio. Our task is to design a portfolio $x$ solving the problem 
\begin{equation}\label{eq:Portfolio}
\min_x \left\{f(x)=-\sum_{t=1}^{T}\ln(r_{t}^{\top}x) : x_{i}\geq 0,\sum_{i=1}^{n}x_{i}=1 \right\}.
\end{equation}
Note that this problem can be cast into an online optimization model \cite{AroHaz06}. It can be seen that $f\in\scrF_{2}$. We remark that this problem appears also in the classical universal prediction problem in information theory and online learning \cite{MerFed98}.

\subsection{The Limits of Standard \ac{FW}}
In this section we explain why the standard analysis of \ac{FW} does not work for \ac{SC} functions. For this, we need to briefly recapitulate the main steps in the modern analysis initiated in \citet{Cla10}, and later on refined by \citet{Haz08} for matrix problems over a spactrahedron, and finally concluded in \citet{Jag13} for general convex compact domains.

The typical merit function employed in \ac{FW} algorithms is the dual gap function 
\begin{equation}\label{eq:gap}
\gap(x)\eqdef\max_{s\in\scrX}\inner{\nabla f(x),x-s}.
\end{equation}
It is easy to see that $\gap(x)\geq 0$ for all $x\in\scrX$, with equality if and only if $x$ is a solution to \eqref{eq:P}.
\begin{algorithm}[t]
 \caption{Standard Frank-Wolfe method}
 \label{alg:FW}
 \begin{algorithmic}
\STATE {\bfseries Input: } $x^{0}\in\dom f\cap \scrX$ initial state, Step size policy $(\alpha_{k})_{k\geq 0}$ (either $\alpha_{k}=\frac{2}{k+2}$, or via line-search); $\eps>0$ tolerance level
\FOR{$k=0,1,\ldots$}
 \IF{$\gap(x^{k})>\eps$}
  \STATE Obtain $s^{k}=s(x^{k})$
  \STATE Set $x^{k+1}=x^{k}+\alpha_{k}(s^{k}-x^{k})$
\ENDIF
\ENDFOR
\end{algorithmic}
\end{algorithm}
The convergence analysis of Algorithm \ref{alg:FW} relies on the \emph{curvature constant} \cite{Jag13}
\begin{align*}
C_{f}\eqdef \sup_{x,s\in\scrX,\gamma\in[0,1]} \frac{2}{\gamma^{2}} D_{f}(x+\gamma(s-x),x),
\end{align*}
where we define the \emph{Bregman divergence} of the smooth convex function $f$ as 
\begin{align*}
D_{f}(y,x)\eqdef f(y)-f(x)-\inner{\nabla f(x),y-x}
\end{align*}
for all $x,y\in\dom f$. Under the assumption that $C_{f}<\infty$, \citet{Jag13} proved sublinear $O(1/k)$ rate of convergence in terms of the dual gap function $\gap(\cdot)$. Unfortunately, minimizing a \ac{SC} function over a compact set does not necessarily give us a finite curvature constant, as the following example illustrates.

\begin{example}\label{example}
Consider the function $f(x,y)=-\ln(x)-\ln(y)$ considered over the set $\scrX=\{(x,y)\in\R^{2}_{+}\vert x+y=1\}$. This function is the standard \ac{SC} barrier for the positive orthant (the log-barrier) and its Bregman divergence is easily calculated as 
\begin{align*}
D_{f}(u,v)=\sum_{i=1}^{2}\left[-\ln\left(\frac{u_{i}}{v_{i}}\right)+\frac{u_{i}}{v_{i}}-1\right]
\end{align*}
where $u=(u_{1},u_{2}),v=(v_{1},v_{2})$. We see that $C_{f}=\infty$. Neither the function $f$, nor its gradient, is Lipschitz continuous over the set of interest. Moreover, if we start from $u^0=(1/4,3/4)$, using the standard step-size $\alpha_k=2/(k+2)$, we have $\alpha_0=1$ and $u^1=s(u^0)=(1,0) \notin \dom f$ and the method fails.
\end{example}
We point out that Example \ref{example} is representative for the class of optimization problems of interest in this work since the logarithm is the prime example for a self-concordant function. It is thus clear that the standard analysis based on finite curvature estimates (or, as a particular case, Lipschitz continuous gradient) cannot be applied to analyze \ac{FW} when applied to \eqref{eq:P}.

\subsection{Basic Estimates}
From \citet{Nes18}, Thm. 5.1.6, we know that if $\dom f$ contains no lines, then $\nabla^{2}f(x)\succ 0$ for all $x\in\dom f$ ($f$ is \emph{non-degenerate}). We shall assume that $f$ is non-degenerate. Define the local norm of $u\in\Rn$ at $x\in\dom f$ as $\norm{u}_{x}\eqdef \sqrt{\inner{\nabla^{2}f(x)u,u}}$. The dual norm at $x\in\dom f$ is defined as $\norm{u}^{\ast}_{x}\eqdef\sqrt{\inner{[\nabla^{2}f(x)]^{-1}u,u}}$. Given $f\in\scrF_{M}$, we define the distance-like function 
$$
\metric(x,y)\eqdef\frac{M}{2}\norm{y-x}_{x},
$$
as well as 
$$
\omega(t)\eqdef t-\ln(1+t),\text{ and }\omega_{\ast}(t)\eqdef-t-\ln(1-t).
$$
It is not hard to verify that $\omega(t)\geq 0$ for all $t>-1$ and $\omega_{\ast}(t)\geq 0$ for every $t<1$; Moreover, $\omega$ and $\omega_{\ast}$ are increasing and strictly convex in $[0,\infty)$ and $[0,1)$, respectively. For all $x,\tilde{x}\in\dom f$ Thms. 5.1.8 and 5.1.9 in \citet{Nes18} state that
\begin{align}
\label{eq:down}
f(\tilde{x})&\geq f(x)+\inner{\nabla f(x),\tilde{x}-x}+\frac{4}{M^{2}}\omega\left(\metric(x,\tilde{x})\right)\\
f(\tilde{x})&\leq f(x)+\inner{\nabla f(x),\tilde{x}-x}+\frac{4}{M^{2}}\omega_{\ast}\left(\metric(x,\tilde{x})\right)
\label{eq:up}
\end{align}
where in the latter inequality we assume $\metric(x,\tilde{x})<1$.\\
An attractive feature of \ac{SC} functions, in particular from the point of view of \ac{FW} algorithms \cite{Jag13}, is that self-concordance is invariant under affine transformations and re-parametrizations of the domain. See Section A of the supplementary materials for a precise statement of this fact, and further properties of \ac{SC} functions.
Affine invariance allows us to cover a very broad class of composite convex optimization problems of the form
\begin{equation}\label{eq:ERM}
\min_{x\in\scrX}\{f(x)\eqdef g(\bE x)+\inner{q,x}\},
\end{equation}
where $g$ is \ac{SC}, and $\bE\in\R^{m\times n},q\in\Rn$. This formulation covers empirical risk-minimization with a convex norm-like regularization, as explained in Section \ref{sec:examples}, among many others.

We have the following existence result for solutions of problem \eqref{eq:P}. For a proof see Section A in the supplementary materials.
\begin{proposition}\label{prop:existence}
Suppose there exists $x\in\dom f\cap\scrX$ such that $\norm{\nabla f(x)}_{x}^{\ast}\leq \frac{2}{M}$. Then \eqref{eq:P} admits a unique solution.
\end{proposition}
We note that there can be at most one solution to problem \eqref{eq:P} because $\nabla^{2}f\succ 0$ on $\dom f\cap\scrX$. Henceforth, we will assume that \eqref{eq:P} admits a solution $x^{\ast}$. 
 
\section{\ac{FW} works for Self-Concordant Functions}
\label{sec:sublinear}
This section gives a high-level overview on the necessary modifications of \ac{FW} in order to derive provably convergent schemes for minimizing an \ac{SC} function over a compact convex set $\scrX$. We propose two new step-size policies, both derived from different approximating local models guaranteeing a version of the celebrated descent lemma. To explain our approach, let us define for all $x\in\dom f$
\begin{equation}\label{eq:e}
\ce(x)\eqdef\metric(x,s(x))=\frac{M}{2}\norm{s(x)-x}_{x}.
\end{equation}
The first step-size rule we analyze, called Variant 1 (V1), is based on the local sufficient decrease property \eqref{eq:up}. In particular, the step-size policy in V1 minimizes at each iteration a local model in the r.h.s. of \eqref{eq:up}, and reads explicitly as 
\begin{equation}
\alpha_{k}^{V1}\eqdef \min\left\{1,\frac{\gap(x^{k})}{\ce(x^{k})(\gap(x^{k})+\frac{4}{M^{2}}\ce(x^{k}))}\right\}.
\end{equation}
Observe, that this step-size rule needs to compute the gap function and the local distance $\ce(x^{k})$. For evaluating this local distance, it is not necessary to evaluate the full Hessian matrix $\nabla^{2}f(x)$, but rather only the matrix-vector product $\nabla^{2}f(x)(s^{k}-x^{k})$. In many situations of interests, this vector can be computed easily. A particular appealing instance where this applies is the generalized linear model $f(x)=\sum_{i=1}^{N}f_{i}(\inner{a_i,x})$, each function $f_{i}$ being \ac{SC}. In this case, the Hessian matrix is a sum of rank-one matrices, Hessian-vector product amounts to vector-vector and vector-scalar products, and the cost of calculating $\ce(x^{k})$ is of the same order as the cost of calculating the gradient.

Variant 2 (V2) is based on a backtracking procedure in the spirit of \citet{PedNegAskJag20}. Specifically, V2 employs the quadratic model 
\begin{equation}\label{eq:Q}
Q(x^{k},t,\mu)\eqdef f(x^{k})-t\gap(x^{k})+\frac{t^{2}\mu}{2}\norm{s(x^{k})-x^{k}}^{2}_{2}
\end{equation}
If $\mu \geq \scrL_{k}$, where $\scrL_{k}$ is a local gradient Lipschitz constant estimate, $Q(x^{k},t,\mu)$ is an upper approximation for the objective. This approximation has the advantage that its minimum over $t\in[0,1]$ can be computed in closed form, yielding the step size
\begin{equation}
\alpha_{k}^{V2}\eqdef \min\left\{1,\frac{\gap(x^{k})}{\scrL_{k}\norm{s(x^{k})-x^{k}}_{2}^{2}}\right\}.
\end{equation}
The quality of this quadratic approximation will depend on the local Lipschitz estimate $\scrL_{k}$. This parameter needs to be carefully selected to ensure convergence while keeping the number of function evaluations to a minimum. This is achieved through Algorithm \ref{alg:backtrack}, which mimics the backtracking procedure of \citet{PedNegAskJag20}.

\begin{algorithm}
 \caption{Adaptive Frank-Wolfe method for \ac{SC} functions}
 \label{alg:AFW} 
 \begin{algorithmic}
\STATE {\bfseries Input: } $x^{0}\in\dom f\cap \scrX$ initial state, $\beta\in(0,1)$, Tolerance $\eps>0$.
\FOR{$k=1,\ldots$}
 \IF{$\gap(x^{k})>\eps$}
  \STATE Obtain $s^{k}=s(x^{k})$ and set $v^{k}=s^{k}-x^{k}$.
  \STATE Obtain $\ce^{k}=\frac{M}{2}\norm{v^{k}}_{x^{k}}$
  \IF{Variant 1:}
  \STATE Set $\ct_{k}=\frac{\gap(x^{k})}{\ce^{k}(\gap(x^{k})+\frac{4}{M^{2}}\ce^{k})}$
  \STATE Set $\alpha_{k}=\min\{1,\ct_{k}\}$
  \ENDIF
  \IF{Variant 2: }
  \STATE Set $(\alpha_{k},\scrL_{k})=\mathtt{step}(f,v^{k},x^{k},\gap(x^{k}),\scrL_{k-1})$ 
  \ENDIF
  \STATE Set $x^{k+1}=x^{k}+\alpha_{k}v^{k}$
\ENDIF
\ENDFOR
\end{algorithmic}
\end{algorithm}
The backtracking performed in this subroutine defines a candidate step size $\alpha$ and checks whether the sufficient decrease condition
\begin{equation}\label{eq:SD}
f(x^{k}+\alpha(s^{k}-x^{k}))<Q(x^{k},\alpha,\mu)
\end{equation}
is satisfied. If not, then we increase the proposed Lipschitz estimate to $\gamma_{u}\scrL,\gamma_{u}>1,$ and repeat. As will be shown in the analysis of this method, this described backtracking procedure is guaranteed to stop in finite steps. 
\begin{algorithm}
 \caption{Function $\mathtt{step}(f,v,x,g,\scrL)$}
 \label{alg:backtrack}
\begin{algorithmic}
    \STATE Choose $\gamma_{u}>1,\gamma_{d}<1$
    \STATE Choose $\mu\in[\gamma_{d}\scrL,\scrL]$
    \STATE $\alpha=\min\{\frac{g}{\scrL\norm{d}^{2}_{2}},1\}$
    \IF{$f(x+\alpha v)>Q(x,\alpha,\mu)$}
   \STATE $\mu=\gamma_{u}\mu$
    \STATE $\alpha=\min\{\frac{g}{\mu\norm{v}^{2}_{2}},1\}$
    \ENDIF
    \STATE Return $\alpha,\mu$
\end{algorithmic}
\end{algorithm}
\subsection{Analysis of Variant 1}
Section B in the supplementary materials contains complete proofs of all the derivations and claims made in this section.

For $x\in\dom f$, define the target vector $s(x)$ as in \eqref{eq:s}, and $\gap(x)$ as in \eqref{eq:gap}. Given $x\in\scrX$ and $t>0$, set $x^{+}_{t}\eqdef x+t(s(x)-x)$. Assume that $\ce(x)\neq 0$. By construction, 
 \begin{align*}
\metric(x,x^{+}_{t})=\frac{tM}{2}\norm{s(x)-x}_{x}=t\ce(x)<1,
 \end{align*}
 iff $t<1/\ce(x)$. Choosing $t\in(0,1/\ce(x))$, we conclude from \eqref{eq:up} 
 \begin{align*}
 f(x^{+}_{t})&\leq f(x)+\inner{\nabla f(x),x^{+}_{t}-x}+\frac{4}{M^{2}}\omega_{\ast}(t\ce(x))\\
 &= f(x)-t\gap(x)+\frac{4}{M^{2}}\omega_{\ast}(t\ce(x))
 \end{align*}
This shows that when minimizing an \ac{SC} function, we can search for a step size $\alpha_{k}$ which minimizes the r.h.s. of the previous inequality by maximizing
 \begin{equation}\label{eq:eta}
 \eta_{x}(t)\eqdef t\gap(x)-\frac{4}{M^{2}}\omega_{\ast}(t\ce(x)),
 \end{equation} 
for  $t\in(0,1/\ce(x)).$ As shown in Section B.2 of the supplementary materials, this function attains a unique maximum at the value 
\begin{equation}\label{eq:t}
 \ct(x)\eqdef \frac{\gap(x)}{\ce(x)(\gap(x)+\frac{4}{M^{2}}\ce(x))}\equiv\frac{\gamma(x)}{\ce(x)}.
 \end{equation}

We then construct the step size sequence $(\alpha_{k})_{k\geq 0}$ by setting $\alpha_{k}=\min\{1,\ct(x^{k})\}$. The convexity of $\scrX$ and the fact that $\alpha_{k}\ce(x^{k})<1$ guarantee that $(x^{k})_{k\geq 0}\subset\dom f\cap\scrX$. Thus, at each iteration, we reduce the objective function value by at least the quantity $\Delta_{k}\equiv\eta_{x^{k}}(\alpha_{k})$, so that $f(x^{k+1})\leq f(x^{k})-\Delta_{k}<f(x^{k}).$ Hence, we see that the method induces a sequence of function values $(f(x^{k}))_{k\geq 0}$ which is monotonically decreasing by at least the amount $\Delta_{k}$ in each iteration. Equivalently,
\[
 (x^{k})_{k\geq 0}\subset\scrS(x^{0})\eqdef\{x\in\dom f\cap\scrX\vert f(x)\leq f(x^{0})\}.
\] 
Furthermore, $\sum_{k\geq 0}\Delta_{k}<\infty$, and hence the sequence $\left(\Delta_{k}\right)_{k\geq 0}$ converges to 0, and 
 \begin{align*}
 \min_{0\leq k<K}\Delta_{k}\leq\frac{1}{K}(f(x^{0})-f^{\ast}).
 \end{align*}
It follows from Lemma B.1 in the supplementary materials that $\frac{4}{M^{2}}\omega(\metric(x^{\ast},x))\leq f(x)-f^{\ast}$ for all $x \in \dom f$. Consequently, 
 \begin{align*}
 \scrS(x^{0})\subseteq\left\{x\in\dom f\vert\omega(\metric(x^{\ast},x))\leq \frac{M^{2}}{4}(f(x^{0})-f^{\ast})\right\}.
 \end{align*}
Since $\omega(\cdot)$ is continuous and increasing, $\scrS(x^{0})$ is a closed, bounded and convex set. Accordingly, defining by $\lambda_{\max}(\nabla^{2}f(x))$ and $\lambda_{\min}(\nabla^{2}f(x))$ the largest, respectively smallest, eigenvalue of the Hessian $\nabla^{2}f(x)$, the numbers 
\begin{align*}
&L_{\nabla f}\eqdef\max_{x\in\scrS(x^{0})}\lambda_{\max}(\nabla^{2}f(x)),\text{ and } \\
&\sigma_{f}\eqdef \min_{x\in\scrS(x^{0})}\lambda_{\min}(\nabla^{2}f(x)),
 \end{align*}
 are well defined and finite. In particular, for all $x\in\scrS(x^{0})$ we have a restricted strong convexity of the function $f$, in the sense that 
 \begin{equation}\label{eq:gbound}
f(x)-f(x^{\ast})\geq \frac{\sigma_{f}}{6}\norm{x-x^{\ast}}_{2}^{2}.
 \end{equation}
 In terms of these quantities, we can bound the sequence $(\ce^{k})_{k\geq 0}$, defined as $\ce^{k}\eqdef \ce(x^{k})$, as
 \begin{equation}\label{eq:normbound}
\frac{M\sqrt{\sigma_{f}}}{2}\norm{s^{k}-x^{k}}_{2}\leq \ce^{k}\leq \frac{M\sqrt{L_{\nabla f}}}{2}\norm{s^{k}-x^{k}}_{2}.
 \end{equation}


In order to derive convergence rates, we need to lower bound the per-iteration decrease in the objective function. A detailed analysis, given in Section B of the supplementary materials, reveals that we can construct a useful minorizing sequence for the per-iteration function decrease $(\Delta_{k})_{k\geq 0}$ as 
\begin{equation}\label{eq:lowDelta}
\Delta_{k}\geq \min\{\ca\gap(x^{k}),\cb\gap(x^{k})^{2}\},
\end{equation}
where $\ca\eqdef \min\left\{\frac{1}{2},\frac{2(1-\ln(2))}{M\sqrt{L_{\nabla f}}\diam(\scrX)}\right\}$ and $\cb\eqdef\frac{1-\ln(2)}{L_{\nabla f}\diam(\scrX)^{2}}$. With the help of this lower bound, we are now able to establish the $\scrO(k^{-1})$ convergence rate in terms of the approximation error $h_{k}\eqdef f(x^{k})-f^{\ast}$. 

By convexity, we have $\gap(x^{k})\geq h_{k}$. Therefore, the lower bound for $\Delta_{k}$ in \eqref{eq:lowDelta} can be estimated as $\Delta_{k}\geq\min\{\ca h_{k},\cb h_{k}^{2}\}$. Hence,
\begin{equation}\label{eq:h}
h_{k+1}\leq h_{k}-\min\{\ca h_{k},\cb h_{k}^{2}\}\qquad\forall k\geq 0.
\end{equation}
Given this recursion, we can identify two phases characterizing the process $(h_{k})_{k\geq 0}$. In Phase I, the approximation error $h_k$ is at least $\ca/\cb$, and in Phase II the approximation error falls below this value. The cut-off value $\ca/\cb$ determines the nature of the recursion \eqref{eq:h}, and yields immediately an estimate for the iteration complexity of Variant 1 of Algorithm \ref{alg:AFW}.
\begin{theorem}\label{thm:1}
For all $\eps>0$, define the stopping time 
\begin{equation}
N_{\eps}(x^{0})\eqdef\inf\{k\geq 0\vert h_{k}\leq\eps\}.
\end{equation}
Then, 
\begin{equation}\label{eq:N}
N_{\eps}(x^{0})\leq \ceil[\bigg]{\frac{1}{\ca}\ln\left(\frac{h_{0}\cb}{\ca}\right)} +\frac{L_{\nabla f}\diam(\scrX)^{2}}{(1-\ln(2))\eps}.
\end{equation}
\end{theorem}
The proof is in Section B.1 of the supplementary materials. 
\subsection{Analysis of Variant 2}
For the analysis of V2, we first need to establish well-posedness of the backtracking scheme $\mathtt{step}(f,v,x,g,\scrL)$. Calling this routine at the position $x=x^{k}$ within the execution of Algorithm \ref{alg:AFW}, we require the search direction $v^{k}\eqdef s^{k}-x^{k}$, where $s^{k}=s(x^{k})$ is the target vector \eqref{eq:s}. Define 
\begin{equation}\label{eq:gamma}
\gamma_{k}\eqdef\max\{t\geq 0\vert x^{k}+t(s^{k}-x^{k})\in\scrS(x^{k})\}
\end{equation}
as the largest step size guaranteeing feasibility and decrease of the objective function value. Hence, for all $t\in[0,\gamma_{k}]$, we have $f(x^{k}+t(s^{k}-x^{k}))\leq f(x^{k})$, and $x^{k}+t(s^{k}-x^{k})\in\scrX$. 

\begin{lemma}
Assume that $x^{k}\in\scrS(x^{0})$ for all $k\geq 0$. Then, for all $t\in[0,\gamma_{k}]$, it holds true that $x^{k}+t(s^{k}-x^{k})\in \scrS(x^{k})$, and 
\[
\norm{\nabla f(x^{k}+t(s^{k}-x^{k}))-\nabla f(x^{k})}_{2}\leq L_{\nabla f}t \norm{s^{k}-x^{k}}_{2}. 
\]
\end{lemma}
Using this estimate, a localized version of the celebrated descent Lemma reads as 
\begin{align*}
f(x^{k}+t(s^{k}-x^{k}))\leq Q(x^{k},t,L_{\nabla f}),
\end{align*}
for all $t\in[0,\gamma_{k}]$ for the quadratic model introduced in \eqref{eq:Q}. A direct minimization strategy for constructing a feasible step size policy based on the majorizing quadratic model $t\mapsto Q(x^{k},t,L_{\nabla f})$ would yield the step size 
\begin{equation}\label{eq:tau}
\tau_{k}(L_{\nabla f})\eqdef\min\left\{1,\frac{\gap(x^{k})}{L_{\nabla f}\norm{s^{k}-x^{k}}^{2}_{2}}\right\}.
\end{equation}
The main problem with this approach is that it needs the global parameter $L_{\nabla f}$. In practice this quantity is often hard to obtain and, furthermore, frequently numerically quite large. This in turn renders the step size to be inevitably small, leading to bad performance of the method. Our solution strategy is thus to implement the adaptive backtracking procedure calling the function $\mathtt{step}(f,v^{k},x^{k},\gap(x^{k}),\scrL_{k-1})$, which requires a local Lipschitz estimate $\scrL_{k-1}$, guaranteed to be smaller than $L_{\nabla f}$. This subroutine needs an initial estimate $\scrL_{-1}$ whose choice is explained in Section C of the supplementary materials.

Furthermore, Algorithm \ref{alg:backtrack} comes with two hyperparameters $\gamma_{d}<1<\gamma_{u}$. Practically efficient values for these parameters are reported to be $\gamma_{d}=0.9$ and $\gamma_{u}=2$ \cite{PedNegAskJag20}. Section C of the supplementary materials gives an upper bound estimate on the number of necessary function evaluations during a single execution of the backtracking procedure.
The $\scrO(k^{-1})$ convergence rate of V2 involves the \emph{dual objective function} 
$$
\psi(z)\eqdef -f^{\ast}(z)-H_{\scrX}(-z),
$$
where $f^{\ast}$ denotes the Fenchel conjugate of $f$ and $H_{\scrX}(c)\eqdef \sup_{x\in\scrX}\inner{x,c}$ is the support function over $\scrX$. Note that $\psi$ is concave and we have $f^{\ast}\eqdef \min_{x\in\scrX}f(x)=\max_{u}\psi(u)$. The following result is then similar to Theorem 3 in \citet{PedNegAskJag20}, and reports the $\scrO(k^{-1})$ convergence rate in terms of the sequence of approximation errors $h_{k}=f(x^{k})-f^{\ast}$.

\begin{theorem}\label{th:V2}
Let $f$ be an \ac{SC} function and $(x^{k})_{k\geq 0}$ generated by Variant 2 of Algorithm \ref{alg:AFW}. 
Then
\begin{equation}
h_{k}\leq\frac{2\gap(x^{0})}{(k+1)(k+2)}+\frac{k\diam(\scrX)^{2}}{(k+1)(k+2)}\bar{\scrL}_{k}
\end{equation}
where $\bar{\scrL}_{k}\eqdef\frac{1}{k}\sum_{i=0}^{k-1}\scrL_{i}$ is the arithmetic mean over all Lipschitz estimates computed by Algorithm \ref{alg:backtrack} during $k$ executions of the main protocol.
\end{theorem}
Since $\scrL_{k}\leq L_{\nabla f}$, this implies $h_{k}=\scrO(k^{-1})$. The proof of this Theorem can be found in supplementary material Section C.3.

\section{Inducing Linear Convergence}
\label{sec:fast}
In this section we use the construction of a \ac{LLOO} given in \citet{GarHaz16} to deduce an accelerated version of the base \ac{FW} Algorithm \ref{alg:AFW}. \citet{GarHaz16} give an explicit construction of such an oracle for \emph{any} polytope of form $\scrX=\{x\in\Rn\vert Ax= a,Bx\leq b\},$ where $A,B\in\R^{m\times n}$ and $a,b\in\R^{m}$. In the particularly important case where $\scrX$ is a standard unit simplex, the \ac{LLOO} can be constructed efficiently. In the supplementary material we describe the construction for this special geometry. 
Let $\ball(x,r)\eqdef\{y\in\Rn\vert \norm{y-x}_{2}\leq r\}$ denote the Euclidean ball with radius $r>0$ and center $x$.

\begin{definition}[\citet{GarHaz16}, Def. 2.5]
A procedure $\scrA(x,r,c)$, where $x\in\scrX,r>0,c\in\Rn,$ is a \ac{LLOO} with parameter $\rho\geq 1$ for the polytope $\scrX$ if $\scrA(x,r,c)$ returns a point $s\in\scrX$ such that for all $y\in\ball(x,r)\cap\scrX$
\begin{equation}
\inner{c,y}\geq\inner{c,s}\text{ and }\norm{x-s}_{2}\leq \rho r.
\end{equation}
\end{definition}

Assuming the availability of a procedure $\scrA(x,r,c)$ for any point $x\in\scrX$, we run Algorithm \ref{alg:LLOO}.
\begin{algorithm}
 \caption{LLOO-based convex optimization}
 \label{alg:LLOO}
 \begin{algorithmic}
\STATE {\bfseries Input: } $\scrA(x,r,c)$-LLOO with parameter $\rho\geq 1$ for polytope $\scrX$, $f\in\scrF_{M}$. $\sigma_{f}>0$ convexity parameter.\\
 $x^{0}\in\dom f\cap \scrX$, and let $h_{0}= f(x^{0})-f^{\ast}$, $c_{0}=1$.
 \FOR{k=0}
\STATE Obtain $r_{0}=\sqrt{\frac{6\gap(x^{0})}{\sigma_{f}}}$.
 \STATE Obtain $s^{0}=\scrA(x^{0},r_{0},\nabla f(x^{0}))$
 \STATE Set $\alpha_{0}=\frac{1}{1+\ce^{0}}\min\{1,\frac{\gap(x^{0})}{\frac{4}{M^{2}}(\ce^{0})^{2}}\}$, where $\ce^{0}=\norm{s^{0}-x^{0}}_{x^{0}}$.
 \STATE Update $x^{1}=x^{0}+\alpha_{0}(s^{0}-x^{0})$ 
 \ENDFOR
\FOR{$k=1,\ldots$}
 \IF{$\gap(x^{k})>\eps$}
\STATE Set $c_{k}=\exp\left(-\frac{1}{2}\sum_{i=0}^{k-1}\alpha_{i}\right)$
\STATE Set $\alpha_{k}=\min\{\frac{c_k\gap(x^0)}{\frac{4}{M^{2}}(\ce^k)^2},1\}\frac{1}{1+\ce^k}$
\STATE Set $r_{k}=r_{0}c_{k}$.
\STATE Obtain $s^{k}=\scrA(x^{k},r_{k},\nabla f(x^{k}))$
\STATE Set $x^{k+1}=x^{k}+\alpha_{k}(s^{k}-x^{k})$
\ENDIF
\ENDFOR
\end{algorithmic}
\end{algorithm}
Our analysis of this Algorithm departs from eq. \eqref{eq:up}, saying that for $\alpha\in(0,1/\ce^{k})$ we have 
\begin{align*}
f(x^{k+1})\leq f(x^{k})+\alpha\inner{\nabla f(x^{k}),s^{k}-x^{k}}+\frac{4}{M^{2}}\omega_{\ast}(\alpha\ce^{k}), 
\end{align*}
To exploit the power of the \ac{LLOO}, we need to control the right-hand side of the previous display by bounding carefully the radius of a single step made by the algorithm. Via a delicate induction argument, based on estimates for \ac{SC} functions and \citet{GarHaz16}, we obtain the announced linear convergence result.
\begin{theorem}
\label{thm:LCON}
Let $(x^{k})_{k\geq 0}$ be generated by Algorithm \ref{alg:LLOO}. Then for all $k\geq 0$ we have $x^{\ast}\in\ball(x^{k},r_{k})$ and 
\begin{equation}\label{eq:fast}
h_{k}\leq \gap(x^{0})\exp\left(-\frac{1}{2}\sum_{i=0}^{k-1}\alpha_{i}\right)
\end{equation}
In particular, defining $\bar{\alpha}\eqdef \min\{ \frac{\sigma_f}{6L_{\nabla f} \rho^2},1\}\frac{1}{1+ \sqrt{L_{\nabla f}}\frac{M\diam(\scrX)}{2}}$, we see that $\alpha_{k}\geq\bar{\alpha}$, and therefore
\[
h_{k}\leq \gap(x^{0})\exp(-k\bar{\alpha}/2).
\]

\end{theorem}
The proof is in Section C of the supplementary materials.

\section{Numerical Experiments}
\label{sec:numerics}
\begin{figure}[t]
\begin{center}
\centerline{\includegraphics[width=1.\columnwidth]{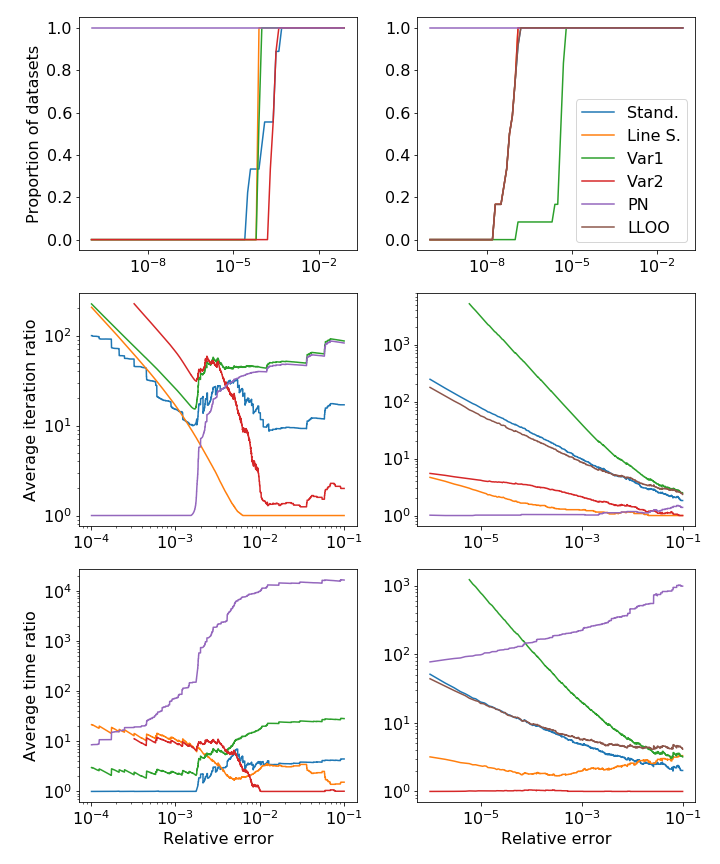}}
\caption{Performance of Algorithm \ref{alg:AFW} and \ref{alg:LLOO} in the Poisson inverse problem \eqref{eq:Poisson} (left column) and in the Portfolio optimization problem \eqref{eq:Portfolio} (right column). The top row reports the fraction of datasets for which a given relative error was achieved by an algorithm. The middle row presents the average ratio between number of iterations taken to reach a certain relative error to the method which reached this error in the minimal number of iterations. The bottom row presents the average time ratio between each method and that of the fastest method in order to reach a certain relative error.}
\label{fig:Portfolio}
\end{center}
\vskip -0.2in
\end{figure}

In the numerical experiments we tested the performance of Variant 1 (V1) and Variant 2 (V2) of Algorithm \ref{alg:AFW}, and compared them with the performance of Frank-Wolfe with {standard} step-size of $\frac{2}{k+2}$ (stand.), and step-size determined by exact {line-search} (Line-S.). As a further benchmark, the self-concordant Proximal-Newton (PN) of \citet{CevKyrTra15}, as implemented in the SCOPT package\footnote{\url{https://www.epfl.ch/labs/lions/technology/scopt/}}, is included. For the portfolio optimization problem, Algorithm \ref{alg:LLOO} is also implemented. All codes are written in Python 3, with packages for scientific computing NumPy 1.18.1 and SciPy 1.4.1. The experiments were conducted on a PC with Intel Core i5-7500 3.4GHzs, with a total of 16GB RAM. In both experiments the Frank-Wolfe based methods have been terminated after 50,000 iterations. Because of its higher computational complexity, we decided to stop PN after 1,000 iterations. Each algorithm was stopped early if the optimality gap in a given iteration was lower than $1e-10$. LLOO was only implemented for the portfolio selection problem, using the local linear oracle given in \cite{GarHaz16} for the simplex, as described in Section F in the supplementary materials.\footnote{The codes are accessible via \url{https://github.com/kamil-safin/SCFW }.}\\
For the portfolio optimization problem we used synthetic data, as in Section 6.4 of \cite{SunTran18}. The details of the data generating process are as follows. We generate matrix $R$ with given price ratios as: $R := ones(n, p) + N (0, 0.1)$, which allows the closing price to vary about 10\% between two consecutive periods. We used different sizes of matrix $R$: $(n, p) = (1000, 800), (1000, 1200)$, and $(1000, 1500)$ with 4 samples for each size. Hence, there are totally 12 datasets. The right column of Figure~\ref{fig:Portfolio} displays the performance of all Algorithms developed in this paper when averaged over all 12 samples. The detailed results for individual samples are reported in Section E of the supplementary materials. 

For the Poisson inverse problem we used the datasets a1a-a9a from the LIBSVM library \cite{LIBSVM}. Averaged performance of the \ac{FW} algorithms over these 9 datasets is displayed in the left column of Figure \ref{fig:Portfolio}. The detailed results for each individual data set are reported in Section E in the supplementary materials.

The comparison between all the algorithms is made by the construction of versions of performance profiles, following \citet{DolMor02}. In order to present the result, we first estimate $f^*$ by the maximal lower bound on the function value achieved by any of the algorithms, and compute the relative error attained by each of the methods at iteration $k$. More precisely, given the set of methods $\scrS$ and test problems $\scrP$, denote by $F_{ij}$ the function value attained by method $i$ on problem $j$. If $(x^{k}_{ij})_{k}$ denotes the sequence produced  by method $i$ on problem $j$, we define the relative error as 
$r^{k}_{ij}=\frac{f(x^{k}_{ij})-\min\{F_{sj}\vert s\in\scrS\}}{\min\{F_{sj}\vert s\in\scrS\}}$.
Now, for all methods $i\in\scrS$ and any relative error $\epsilon$, The  top row of Figure~\ref{fig:Portfolio} shows the proportion of datasets that achieves a relative error of $\epsilon$, that is $\rho_i(\epsilon)= \frac{1}{|\scrP|}\abs{\{j\in\scrP\vert \exists k,\; r^{k}_{ij}\leq \eps\}}$.
We are also interested in comparing iteration complexity and CPU time between the methods. For that purpose, we define 
 $N_{ij}(\eps)=\min\{k\geq 0\vert r^{k}_{ij}\leq \eps\}$ and $T_{ij}(\eps)$ as the first iteration and the minimal CPU time in which method $i\in\scrS$ to achieves a relative error $\eps$ on problem $j\in\scrP$. 
In the second row of Figure~\ref{fig:Portfolio}, we use the following average iteration ratio $\tilde{\rho}_{i}(\eps)=\frac{1}{|\scrP|}\sum_{j\in\scrP}\frac{N_{ij}(\eps)}{\min\{N_{sj}(\eps)\vert s\in\scrS\}}$ for comparing the iteration complexity of all the methods. 
In the third row of Figure~\ref{fig:Portfolio}, we display the average time ratio $\hat{\rho}_{i}(\eps)=\frac{1}{|\scrP|}\sum_{j\in\scrP}\frac{T_{ij}(\eps)}{\min\{T_{sj}(\eps)\vert s\in\scrS\}}$ for comparing the computational time of all the methods. In both cases, as the average ratio is closer to 1 the performance of the method is closer to the best performance.

As expected, the top row of Figure~\ref{fig:Portfolio} show that PN obtains a lower relative error than any of the Frank-Wolfe based methods. However, for the Frank-Wolfe methods, none of the step-sizes choices has a clear advantage over the other in obtaining lower relative error values in both examples. Moreover, in the portfolio, while PN has the lowest iteration complexity, Variant 2 in fact achieves the best times for all values of relative error, followed by line-search, with LLOO matching the performance of the standard step-size. We remark that the standard step-size policy has no theoretical convergence guarantees, as the problem has no finite curvature constant (nor has a Lipschitz continuous gradient over the unit simplex).\\
In the Poisson inverse problem displayed in the left column of Figure \ref{fig:Portfolio}, Variant 2 obtains the best times for higher values of relative errors. The standard $2/(k+2)$ step size rule seems to be a good alternative here. We also implemented the convergent version with the step size $2/(k+3)$ derived in \citet{Odor16} and did not observe any significant difference in practice with the standard step-size policy. It is remarkable that the worst-case step-size policy V1 is competitive with the problem-specific $2/(k+3)$ policy.

\section{Conclusion}
\ac{FW} is a much appraised \ac{FOM} for minimizing smooth convex functions over convex compact sets. The main merit of this method is the relative simplicity of its implementation and projection-free iterations. This yields great scalability properties making it a very attractive method for large-scale optimization. The price to pay for iteration-simplicity are, in general, slow convergence rates, and some sort of bounded curvature assumption.
In this work, we show that for \ac{SC} functions, which are neither strongly convex, nor have bounded curvature, we can obtain a novel step-size policy which, when coupled with local linear minimization oracles, features linear convergence rates. Under more standard assumptions of the feasible set given by \ac{LMO}, we provide two novel step-size policies which lead to standard sublinear convergence rate for minimizing general \ac{SC} functions by \ac{FW} algorithm.
In the future we plan to extend the results to the class of generalized self-concordant functions, as recently defined in \cite{SunTran18}. Other directions of interest for future research are inertial methods, and stochastic optimization. We will address these problems in the near future. 
An interesting avenue could be to incorporate away steps. This seems to be challenging since vanilla away step is not trivially to implement as the function might not be bounded and the away step may lead to a point outside of the domain of the objective.

\section*{Acknowledgements}
The authors sincerely thank Shoham Sabach for his contribution in the early stages of this project, including his part in developing the basic ideas used in this paper. We would also like to thank Quoc Tran-Dinh for sharing MATLAB codes. This research is supported by the COST Action CA16228 "European Network for Game Theory".

\appendix
\section*{Outline} 
The supplementary material of this paper is organized as follows. 
\begin{itemize}
\item Appendix \ref{sec:prop} contains further details on \ac{SC} functions.
\item Appendix \ref{sec:V1} is devoted to proof of Theorem 3.1  for Variant 1 of Algorithm 2. Since this proof relies on some standard estimates on self-concordant functions, we include those auxiliary estimates as well. 
\item Appendix \ref{sec:V2} is organized around the convergence proof of Variant 2 of Algorithm 2, which is Theorem 3.4 in the main text. We also give some guidelines how the parameters and initial values of the backtracking subroutine are chosen. 
\item Appendix \ref{sec:LLOO} contains the linear convergence proof under the availability of the restricted local linear minimization oracle (LLOO)
\item Appendix \ref{sec:numerics} collects detailed evaluations of the numerical performances of the Algorithms constructed in this paper in the context of the Portfolio optimization and the Poisson inverse problem. 
\item Appendix \ref{sec:Simplex} outlines the construction of the LLOO for simplex constraints, following \cite{GarHaz16}.
\end{itemize}

\section{Proofs of Section 2}
\label{sec:prop}

We first introduce a classical result on SC functions, showing its affine invariance. 

\begin{lemma}[\citet{Nes18}, Thm. 5.1.2] 
\label{lem:invariant}
Let $f\in\scrF_{M}$ and $\scrA(x)=Ax+b:\Rn\to\R^{p}$ a linear operator. Then $\tilde{f}\eqdef f\circ\scrA\in\scrF_{M}$.
\end{lemma}
When we apply \ac{FW} to the minimization of a function $f\in\scrF_{M}$, the search direction at position $x$ is determined by the target state $s(x)=s$ defined in \eqref{eq:s}. If $A:\tilde{\scrX}\to\scrX$ is a surjective linear re-parametrization of the domain $\scrX$, then the new optimization problem $\min_{\tilde{\scrX}}\tilde{f}(\tilde{x})=f(A\tilde{x})$ is still within the frame of problem (P). Furthermore, the updates produced by \ac{FW} are not affected by this re-parametrization since 
$$
\inner{\nabla\tilde{f}(\tilde{x}),\hat{s}}=\inner{\nabla f(A\tilde{x}),A\hat{s}}=\inner{\nabla f(x),s}
$$
for $x=A\tilde{x}\in\scrX,s=A\hat{s}\in\scrX$.

\begin{proposition}\label{prop:existence}
Suppose there exists $x\in\dom f\cap\scrX$ such that $\norm{\nabla f(x)}_{x}^{\ast}\leq \frac{2}{M}$. Then (P) admits a unique solution.
\end{proposition}
\begin{proof}
For all $x,y\in\dom f$ we know that 
\begin{align*}
f(y)&\geq f(x)+\inner{\nabla f(x),y-x}+\frac{4}{M^{2}}\omega\left(\frac{M}{2}\norm{y-x}^{2}_{x}\right)\\
&\leq f(x)-\norm{\nabla f(x)}^{\ast}_{x}\cdot\norm{y-x}_{x}\\
&+\frac{4}{M^{2}} \omega\left(\frac{M}{2}\norm{y-x}^{2}_{x}\right)\\
&=f(x)+\left(\frac{2}{M}-\norm{\nabla f(x)}^{\ast}_{x}\right)\norm{y-x}_{x}\\
&-\frac{4}{M^{2}}\ln\left(1+\frac{M}{2}\norm{y-x}_{x}\right).
\end{align*}
Define the level set $\scrL_{f}(\alpha)\eqdef\{x\vert f(x)\leq \alpha\}$ and pick $y\in\scrL_{f}(f(x))$. For such a point, we get 
\begin{align*}
\frac{4}{M^{2}}&\ln\left(1+\frac{M}{2}\norm{y-x}_{x}\right)\\
&\geq\left(\frac{2}{M}-\norm{\nabla f(x)}^{\ast}_{x}\right)\norm{y-x}_{x}.
\end{align*}
Consider the function $t\mapsto \varphi(t)\eqdef\frac{\ln(1+t)}{t}$ for $t>0$. For $t>0$ it is true that $\varphi(t)<1$, and so we need that $\norm{\nabla f(x)}^{\ast}_{x}\leq\frac{2}{M}$. Since $\lim_{t\to\infty}\varphi(t)=0$, it follows that $\scrL_{f}(f(x))$ is bounded. By the Weierstrass theorem, existence of a solution follows (see e.g. \cite{Ber99}). If $x^{\ast}\in\dom f\cap\scrX$ is a solution, we know that 
\[
f(x)\geq f(x^{\ast})+\frac{4}{M^{2}}\omega(\frac{M}{2}\norm{x-x^{\ast}}_{x^{\ast}}).
\]
Hence, if $x$ would be any alternative solution, we immediately conclude that $x=x^{\ast}$.
\end{proof}

\section{Proofs of convergence of Variant 1 of Algorithm 2}
\label{sec:V1}
This supplementary material contains all results needed to establish the convergence of Version 1 of Algorithm 2. We start with some basic estimates helping to proof the main result about this numerical scheme. 
\subsection{Preliminary Results}
\label{sec:prelims}
We recall the basic inequalities for \ac{SC} functions. 
\begin{align}
\label{eq:down}
f(\tilde{x})&\geq f(x)+\inner{\nabla f(x),\tilde{x}-x}+\frac{4}{M^{2}}\omega\left(\metric(x,\tilde{x})\right)\\
f(\tilde{x})&\leq f(x)+\inner{\nabla f(x),\tilde{x}-x}+\frac{4}{M^{2}}\omega_{\ast}\left(\metric(x,\tilde{x})\right)
\label{eq:up}
\end{align}
We need a preliminary error bound around the unique solution.
\begin{lemma}\label{lem:w1}
For all $x\in\dom f$ we have:
\begin{align*}
\frac{4}{M^{2}}\omega\left(\frac{M}{2}\norm{x-x^{\ast}}_{x^{\ast}}\right)\leq f(x)-f(x^{\ast}). 
\end{align*}
\end{lemma}
\begin{proof}
If $x\in\dom f\cap\scrX$, then eq. \eqref{eq:down} shows
\begin{align*}
f(x)&\geq f(x^{\ast})+\inner{\nabla f(x^{\ast}),x-x^{\ast}}+\frac{4}{M^{2}}\omega\left(\metric(x^{\ast},x)\right)\\
&\geq f(x^{\ast})+\frac{4}{M^{2}}\omega(\metric(x^{\ast},x)).
\end{align*}
\end{proof} 

We next pro\pd{ve} a restricted strong convexity property of SC functions. 
 \begin{lemma}\label{lem:gbound}
For all $x\in\scrS(x^{0})$ we have 
 \begin{equation}\label{eq:gbound}
f(x)-f(x^{\ast})\geq \frac{\sigma_{f}}{6}\norm{x-x^{\ast}}_{2}^{2}.
 \end{equation}
 \end{lemma}
\begin{proof}
Lemma \ref{lem:w1} gives $f(x)-f(x^{\ast})\geq\frac{4}{M^{2}}\omega(\metric(x^{\ast},x))$. Observe that for all $t\in[0,1]$
 \begin{align*}
\omega(t)=t-\ln(1+t)=\sum_{j=2}^{\infty}\frac{(-1)^{j}t^{j}}{j}\geq \frac{t^{2}}{2}-\frac{t^{3}}{3}\geq\frac{t^{2}}{6}.
 \end{align*}
Coupled with the fact that $x^{\ast}\in\scrS(x^{0})$ and the hypothesis that $x\in\scrS(x^{0})$, we see that 
\[
f(x)-f^{\ast}\geq \frac{1}{6}\norm{x-x^{\ast}}^{2}_{x^{\ast}}\geq\frac{\sigma_{f}}{6}\norm{x-x^{\ast}}^{2}_{2}.
\]
\end{proof}
Also, we need the next classical fact for SC functions.

\begin{lemma}\label{lem:Dikin}
Let $\scrW(x,r)=\{x'\in\Rn\vert\frac{M}{2}\norm{x'-x}_{x}<r\}$ denote the Dikin ellipsoid with radius $r$ around $x$. For all $x\in \dom f$ we have $\scrW(x,1)\subset\dom f$.
\end{lemma}
\begin{proof}
See \citet{Nes18}.
\end{proof}

\subsection{Estimates for the Algorithm}
For $x\in\dom f$, define the target vector 
\begin{equation}\label{eq:s}
s(x)=\argmin_{x\in\scrX}\inner{\nabla f(x),x},
\end{equation}
and 
\begin{equation}\label{eq:gap}
\gap(x)=\inner{\nabla f(x),x-s(x)}.
\end{equation}
Moreover, for all $x\in\dom f$, let us define 
\begin{equation}\label{eq:e}
\ce(x)\eqdef\metric(x,s(x))=\frac{M}{2}\norm{s(x)-x}_{x}.
\end{equation}
 Given $x\in\scrX$ and $t>0$, set $x^{+}_{t}\eqdef x+t(s(x)-x)$. Assume that $\ce(x)\neq 0$. By construction, 
 \begin{align*}
\metric(x,x^{+}_{t})=\frac{tM}{2}\norm{s(x)-x}_{x}=t\ce(x)<1,
 \end{align*}
 iff $t<1/\ce(x)$. Choosing $t\in(0,1/\ce(x))$, we conclude from \eqref{eq:up} 
 \begin{align*}
 f(x^{+}_{t})&\leq f(x)+\inner{\nabla f(x),x^{+}-x}+\frac{4}{M^{2}}\omega_{\ast}(t\ce(x))\\
 &\leq f(x)-t\gap(x)+\frac{4}{M^{2}}\omega_{\ast}(t\ce(x))
 \end{align*}
This reveals the interesting observation that for minimizing an SC-function, we can search for a step size $\alpha_{k}$ which minimizes the model function
 \begin{equation}\label{eq:eta}
 \eta_{x}(t)\eqdef t\gap(x)-\frac{4}{M^{2}}\omega_{\ast}(t\ce(x)),
 \end{equation} 
defined for $t\in(0,1/\ce(x))$
 \begin{proposition}\label{prop:descent}
 For $x\in\dom f\cap\scrX$, the function $t\mapsto \eta_{x}(t)$ defined in \eqref{eq:eta} is concave and uniquely maximized at the value 
\begin{equation}\label{eq:t}
 \ct(x)\eqdef \frac{\gap(x)}{\ce(x)(\gap(x)+\frac{4}{M^{2}}\ce(x))}\equiv\frac{\gamma(x)}{\ce(x)}.
 \end{equation}
 If $\alpha(x)\eqdef\min\{1,\ct(x)\}$  is used as a step-size in Variant 1 of Algorithm 2, and define $\Delta(x)\eqdef\eta_{x}(\alpha(x))$, then 
 \begin{equation}\label{eq:descent1}
 f(x+\alpha(x)(s(x)-x))\leq f(x)-\Delta(x).
 \end{equation}
\end{proposition}
\begin{proof}
For $x\in\dom f\cap\scrX$, define 
\begin{equation}
\eta_{x}(t)\eqdef t\gap(x)-\frac{4}{M^{2}}\omega_{\ast}(t\ce(x)).
\end{equation} 
We easily compute $\eta''_{x}(t)=\frac{4}{M^{2}}\frac{\ce(x)}{(1-t\ce(x))^{2}}>0$. Hence, the function is concave and uniquely maximized at 
\begin{equation}
\ct(x)\eqdef \frac{\gap(x)}{\ce(x)(\gap(x)+\frac{4}{M^{2}}\ce(x))}\equiv\frac{\gamma(x)}{\ce(x)}.
\end{equation}
Furthermore, one can easily check that $\eta_{x}(0)=0$, and $\eta_{x}(\ct(x))=\frac{4}{M^{2}}\omega\left(\frac{M^{2}}{4}\frac{\gap(x)}{\ce(x)}\right)>0$, whenever $\ce(x)>0$. Hence, it follows that 
\begin{equation}\label{eq:eta-p}
\eta_{x}(t)>0\quad\forall t\in(0,\ct(x)]. 
\end{equation}
\end{proof}
We now construct the step size sequence $(\alpha_{k})_{k\geq 0}$ by setting $\alpha_{k}=\min\{1,\ct(x^{k})\}$ for all $k\geq 0$. Convexity of $\scrX$ and the fact that $\alpha_{k}\ce(x^{k})<1$ guarantees that $(x^{k})_{k\geq 0}\subset\dom f\cap\scrX$. For the feasibility, we use Lemma \ref{lem:Dikin}. Thus, at each iteration, we reduce the objective function value by at least the quantity $\Delta_{k}\equiv\eta_{x^{k}}(\alpha_{k})$, so that $f(x^{k+1})\leq f(x^{k})-\Delta_{k}<f(x^{k}).$ 
\begin{proposition}
\label{prop:Delta}
The following assertions hold for Variant 1 of Algorithm 2:
\begin{itemize}
\item[(a)] $\left(f(x^{k})\right)_{k\geq 0}$ is non-increasing;
\item[(b)] $\sum_{k\geq 0}\Delta_{k}<\infty$, and hence the sequence $\left(\Delta_{k}\right)_{k\geq 0}$ converges to 0;
\item[(c)] For all $K\geq 1$ we have $\min_{0\leq k<K}\Delta_{k}\leq\frac{1}{K}(f(x^{0})-f^{\ast})$.
\end{itemize}
\end{proposition}
\begin{proof}
This proposition can be deduced from Proposition 5.2 in \cite{HBA-GSC}. We give a proof for the sake of being self-contained. Evaluating eq. (9) along the iterate sequence, and calling $\Delta_{k}=\Delta(x^{k})$, we get for all $k\geq 0$,
\begin{align*}
f(x^{k+1})-f(x^{k})\leq -\Delta_{k}.
\end{align*}
Telescoping this expression shows that for all $K\geq 1$, 
\begin{align*}
f(x^{K})-f(x^{0})\leq -\sum_{k=0}^{K-1}\Delta_{k}. 
\end{align*}
Since $\Delta_{k}>0$ , the sequence $\left(f(x^{k})\right)_{k\geq 0}$ is monotonically decreasing. We conclude that for all $K\geq 1$, 
\begin{equation}
\sum_{k=0}^{K-1}\Delta_{k}\leq f(x^{0})-f(x^{K})\leq f(x^{0})-f^{\ast}
\end{equation}
and therefore,
\begin{equation}
\min_{1\leq k\leq K}\Delta_{k}\leq\frac{1}{K}(f(x^{0})-f^{\ast}).
\end{equation}
Hence, $\lim_{k\to\infty}\Delta_{k}=0$. 
\end{proof}

We can bound the sequence $(\ce^{k})_{k\geq 0}$, defined as $\ce^{k}\eqdef \ce(x^{k})$, as
 \begin{equation}\label{eq:normbound}
\frac{M\sqrt{\sigma_{f}}}{2}\norm{s^{k}-x^{k}}_{2}\leq \ce^{k}\leq \frac{M\sqrt{L_{\nabla f}}}{2}\norm{s^{k}-x^{k}}_{2}.
 \end{equation}


In order to derive convergence rates, we need to lower bound the per-iteration decrease in the objective function. A detailed analysis of the sequence $(\Delta_{k})_{k\geq 0}$ reveals an explicit lower bound on the per-iteration decrease which relates the gap function to the sequence $(\Delta_{k})_{k\geq 0}$.
\begin{lemma}\label{lem:delta}
For all $k\geq 0$ we have 
\begin{equation}
\Delta_{k}\geq \min\{\ca\gap(x^{k}),\cb\gap(x^{k})^{2}\},
\end{equation}
where $\ca\eqdef \min\left\{\frac{1}{2},\frac{2(1-\ln(2))}{M\sqrt{L_{\nabla f}}\diam(\scrX)}\right\}$ and $\cb\eqdef\frac{1-\ln(2)}{L_{\nabla f}\diam(\scrX)^{2}}$.
\end{lemma}
\begin{proof}
Let us start with an iteration $k$ at which $\alpha_{k}=\ct(x^{k})$. In this case, we make progress according to 
\begin{align*}
\eta_{x^{k}}(\ct(x^{k}))&=\frac{\gap(x^{k})}{\ce(x^{k})}\gamma(x^{k})+\frac{4}{M^{2}}\gamma(x^{k})\\
&+\frac{4}{M^{2}}\ln\left(\frac{(4/M^{2})\ce(x^{k})}{\gap(x^{k})+(4/M^{2})\ce(x^{k})}\right).
\end{align*}
Define $y:=\frac{(4/M^{2})\ce(x^{k})}{\gap(x^{k})}$. Rewriting the above display in terms of this new variable, we arrive, after some elementary algebra, at the expression
\begin{align*}
\eta_{x^{k}}(t(x^{k}))=\frac{\gap(x^{k})}{\ce(x^{k})}\left[1+y\ln\left(\frac{y}{1+y}\right)\right].
\end{align*}
Consider the function $\phi:(0,\infty)\to(0,\infty)$, given by  $\phi(t)\eqdef 1+ t \ln\left(\frac{t}{1+t}\right)$. When $t \in (0,1)$, since 
\begin{align*}
\phi'(t) &= \ln\left(\frac{t}{1+t}\right) + t \frac{1+t}{t}\left(\frac{1}{1+t}-\frac{t}{(1+t)^2}  \right) \\
& = \ln\left(\frac{t}{1+t}\right) + 1 - \frac{t}{1+t}\\
&= \ln\left(1- \frac{1}{1+t}\right)  + \frac{1}{1+t} < 0,
\end{align*}
we conclude that $\phi(t)$ is decreasing for $t\in(0,1)$. Hence, $\phi(t) \geq \phi(1) = 1-\ln 2$, for all $t \in (0,1)$. On the other hand, if $t \geq 1$, 
\begin{align*}
\frac{\dif}{\dif t}\left(\frac{\phi(t)} {1/t}\right)&=\frac{\dif}{\dif t}(t\phi(t) )\\
&=1 + 2t \ln\left(\frac{t}{1+t}\right)  + \frac{t}{1+t}\geq 0.
\end{align*}
Hence, $t\mapsto \frac{\phi(t)} {1/t}$ is an increasing function for $t\geq 1$, and thus $\phi(t) \geq \frac{1-\ln 2}{t}$, for all $t\geq 1$. We conclude that 
\begin{align*}
\eta_{x^{k}}(\ct(x^{k}))\geq \frac{\gap(x^{k})}{\ce(x^{k})}(1-\ln(2))\min\left\{1,\frac{\gap(x^{k})}{(4/M^{2})\ce(x^{k})}\right\}. 
\end{align*}
Now consider an iteration $k$ in which $\alpha_{k}=1$. The per-iteration decrease of the objective function is explicitly given by 
\begin{align*}
\eta_{x^{k}}(1)&=\left[\gap(x^{k})+\frac{4}{M^{2}}\ce(x^{k})\right]+\frac{4}{M^{2}}\ln(1-\ce(x^{k}))\\
&=\gap(x^{k})\left[1+y+\frac{y}{\ce(x^{k})}\ln(1-\ce(x^{k}))\right].
\end{align*}
Since $\alpha_{k}=1$, it is true that $\ce(x^{k})<\gamma(x^{k})<1$, and therefore $\frac{1}{\ce(x^{k})}\ln(1-\ce(x^{k}))>\frac{1}{\gamma(x^{k})}\ln(1-\gamma(x^{k}))$. Finally, using the identity $1+y=\frac{1}{\gamma(x^{k})}$, we arrive at the lower bound
\begin{align*}
\eta_{x^{k}}(1)&\geq \gap(x^{k})\left[1+y+y(1+y)\ln\left(\frac{y}{1+y}\right)\right]\\
&\geq\frac{\gap(x^{k})}{2}. 
\end{align*}
Summarizing all these computations, we see that for all $k\geq 0$, the per-iteration decrease is at least 
\begin{align*}
&\Delta_{k}\geq\\
&\min\left\{\frac{\gap(x^{k})}{2},\frac{(1-\ln(2))\gap(x^{k})}{\ce(x^{k})},\frac{(1-\ln(2))\gap(x^{k})^{2}}{(4/M^{2})\ce(x^{k})^{2}}\right\}.
\end{align*}
From eq. \eqref{eq:normbound}, we deduce that $\ce(x)\leq\frac{M\sqrt{L_{\nabla f}}}{2}\diam(\scrX)$. Hence, after setting $\ca\eqdef \min\left\{\frac{1}{2},\frac{2(1-\ln(2))}{M\sqrt{L_{\nabla f}}\diam(\scrX)}\right\}$ and $\cb\eqdef\frac{1-\ln(2)}{L_{\nabla f}\diam(\scrX)^{2}}$, we see that 
\[
\Delta_{k}\geq \min\{\ca\gap(x^{k}),\cb\gap(x^{k})^{2}\}.
\]
\end{proof}

\subsection{Proof of Theorem 3.6}
With the help of the lower bound in Lemma \ref{lem:delta}, we are now able to establish the $\scrO(k^{-1})$ convergence rate in terms of the approximation error $h_{k}\eqdef f(x^{k})-f^{\ast}$. 

By convexity, we have $\gap(x^{k})\geq h_{k}$. Therefore, the lower bound for $\Delta_{k}$ can be estimated as $\Delta_{k}\geq\min\{\ca h_{k},\cb h_{k}^{2}\}$, which implies
\begin{equation}\label{eq:h}
h_{k+1}\leq h_{k}-\min\{\ca h_{k},\cb h_{k}^{2}\}\qquad\forall k\geq 0.
\end{equation}
Given this recursion, we can identify two phases characterizing the process $(h_{k})_{k\geq 0}$. In Phase I, the approximation error is at least $\ca/\cb$, and in Phase II the approximation error falls below this value.

For fixed initial condition $x^{0}\in\dom f\cap\scrX$, we can subdivide the time domains according to Phases I and II as
\begin{align*}
&\scrK_{1}(x^{0})\eqdef \{k\geq 0\vert h_{k}>\frac{\ca}{\cb}\},\text{ (Phase I) }\\
&\scrK_{2}(x^{0})\eqdef\{k\geq 0\vert h_{k}\leq \frac{\ca}{\cb}\},\text{ (Phase II)}.
\end{align*}
Since $(h_{k})_{k}$ is monotonically decreasing and bounded from below by the positive constant $\ca/\cb$ on Phase I, the set $\scrK_{1}(x^{0})$ is at most finite. Let us set 
\begin{equation}
T_{1}(x^{0})\eqdef\inf\{k\geq 0\vert h_{k}\leq\frac{\ca}{\cb}\},
\end{equation}
the first time at which the process $(h_{k})$ enters Phase II. To get a worst-case estimate on this quantity, assume that $0\in \scrK_{1}(x^{0})$, so that $\scrK_{1}(x^{0})=\{0,1,\ldots,T_{1}(x^{0})-1\}$. Then, for all $k=1,\ldots,T_{1}(x^{0})-1$ we have $\frac{\ca}{\cb}<h_{k}\leq h_{k-1}-\min\{\ca h_{k-1},\cb h_{k-1}^{2}\}=h_{k-1}-\ca h_{k-1}$. Note that $\ca\leq 1/2$, so we make progressions like a geometric series. Hence, 
$h_{k}\leq (1-\ca)^{k}h_{0}$ for all $k=0,\ldots,T_{1}(x^{0})-1$. By definition $h_{T_{1}(x^{0})-1}>\frac{\ca}{\cb}$, so we get 
$\frac{\ca}{\cb}\leq h_{0} (1-\ca)^{T_{1}(x^{0})-1}$ iff $(T_{1}(x^{0})-1)\ln(1-\ca)\geq \ln\left(\frac{\ca}{h_{0}\cb}\right)$. Hence, 
\begin{equation}\label{eq:T1}
T_{1}(x^{0})\leq \ceil[\bigg]{\frac{\ln\left(\frac{\ca}{h_{0}\cb}\right)}{\ln(1-\ca)}} \leq \ceil[\bigg]{\frac{1}{\ca}\ln\left(\frac{h_{0}\cb}{\ca}\right)}.
\end{equation}
After these number of iterations, the process will enter Phase II, at which $h_{k}\leq\frac{\ca}{\cb}$ holds. Therefore, $h_{k}\geq h_{k+1}+\cb h_{k}^{2}$, or equivalently, 
\begin{align*}
\frac{1}{h_{k+1}}\geq \frac{1}{h_{k}}+\cb\frac{h_{k}}{h_{k+1}}\geq \frac{1}{h_{k}}+\cb.
\end{align*}
Pick $N>T_{1}(x^{0})$ an arbitrary integer. Summing this relation for $k=T_{1}(x^{0})$ up to $k=N-1$, we arrive at 
\begin{align*}
\frac{1}{h_{N}}\geq \frac{1}{h_{T_{1}(x^{0})}}+\cb(N-T_{1}(x^{0})+1).
\end{align*}
By definition $h_{T_{1}(x^{0})}\leq \frac{\ca}{\cb}$, so that for all $N>T_{1}(x^{0})$, we see
\begin{align*}
\frac{1}{h_{N}}\geq \frac{\cb}{\ca}+\cb(N-T_{1}(x^{0})+1).
\end{align*}
Consequently, 
\begin{align*}
h_{N}&\leq \frac{1}{\frac{\cb}{\ca}+\cb(N-T_{1}(x^{0})+1)}\\
&\leq\frac{1}{\cb(N-T_{1}(x^{0})+1)}\\
&=\frac{L_{\nabla f}\diam(\scrX)^{2}}{(1-\ln(2))(N-T_{1}(x^{0})+1)}.
\end{align*}
Define the stopping time $N_{\eps}(x^{0})\eqdef\inf\{k\geq 0\vert h_{k}\leq\eps\}.$ Then, by definition, it is true that $h_{N_{\eps}(x^{0})-1}>\eps$, and consequently, evaluating the bound for $h_{N}$ at $N=N_{\eps}(x^{0})-1$, we obtain the relation
\begin{align*}
\eps\leq \frac{L_{\nabla f}\diam(\scrX)^{2}}{(1-\ln(2))(N_{\eps}(x^{0})-T_{1}(x^{0}))}.
\end{align*}
Combining with the estimate \eqref{eq:T1}, and solving the previous relation of $N_{\eps}(x^{0})$ gives us
\begin{equation}\label{eq:N}
N_{\eps}(x^{0})\leq \ceil[\bigg]{\frac{1}{\ca}\ln\left(\frac{h_{0}\cb}{\ca}\right)} +\frac{L_{\nabla f}\diam(\scrX)^{2}}{(1\pd{-}\ln(2))\eps}.
\end{equation}

\section{Proofs for Variant 2 of Algorithm 2}
\label{sec:V2}
In this section we describe them main steps in the convergence analysis of Variant 2 of Algorithm 2. In order to ensure that the evaluation of the function $\mathtt{step}(f,v,x,g,\scrL)$ needs only finitely many iterations, we need to establish a conceptual global descent lemma. Such a descent property is established in the next Lemma, which corresponds to Lemma 3.2 in the main text. 

\begin{lemma}
Assume that $x^{k}\in\scrS(x^{0})$ for all $k\geq 0$. For all $t\in[0,\gamma_{k}]$, it holds true that $x^{k}+t(s^{k}-x^{k})\in \scrS(x^{k})$, and 
\[
\norm{\nabla f(x^{k}+t(s^{k}-x^{k}))-\nabla f(x^{k})}\leq L_{\nabla f}t \norm{s^{k}-x^{k}}_{2}. 
\]
\end{lemma}
\begin{proof}
The descent property $x^{k}+t(s^{k}-x^{k})\in\scrS(x^{k})$ for $t\in[0,\gamma_{k}]$ follows directly from the definition of $\gamma_{k}$. By the mean-value theorem, for all $\sigma>0$ such that $x^{k}+t(s^{k}-x^{k})\in\scrS(x^{k})$, we have 
\begin{align*}
&\norm{\nabla f(x^{k}+t (s^{k}-x^{k}))-\nabla f(x^{k})}_{2}\\
&=\norm{\int_{0}^{t}\nabla^{2}f(x^{k}+\tau (s^{k}-x^{k}))\dif\tau \cdot(s^{k}-x^{k})}_{2}\\
&\leq \int_{0}^{t}\norm{\nabla^{2}f(x^{k}+\tau (s^{k}-x^{k}))(s^{k}-x^{k})}_{2}\dif\tau\\
&\leq L_{\nabla f}t \norm{s^{k}-x^{k}}_{2}.
\end{align*}
\end{proof}
This implies a localized version of the descent Lemma, which reads as 
\begin{equation}\label{eq:SD}
\begin{array}{l}
f(x^{k}+t(s^{k}-x^{k}))-f(x^{k})\\
-\inner{\nabla f(x^{k}),t(s^{k}-x^{k})}\leq \frac{L_{\nabla f}t^{2}}{2}\norm{s^{k}-x^{k}}_{2}
\end{array}
\end{equation}
for all $t\in[0,\gamma_{k}].$ Introducing the quadratic model
\begin{equation}\label{eq:Q}
Q(x^{k},t,\mu)\eqdef f(x^{k})-t\gap(x^{k})+\frac{t^{2}\mu}{2}\norm{s(x^{k})-x^{k}}^{2}_{2},
\end{equation}
this reads as 
\begin{equation}
f(x^{k}+t(s^{k}-x^{k}))\leq Q(x^{k},t,L_{\nabla f}).
\end{equation}
\subsection{Initial parameters}
The backtracking subroutine, Algorithm 3, needs to know initial values for the Lipschitz estimate $\scrL_{-1}$. In \citet{PedNegAskJag20}, it is recommended to use the following heuristic: Choose $\eps=10^{-3}$, or any other positive numbers smaller than this. Then set 
\[
\scrL_{-1}=\frac{\norm{\nabla f(x^{0})-\nabla f(x^{0}+\eps(s^{0}-x^{0}))}}{\eps\norm{s^{0}-x^{0}}}
\]
The function $\mathtt{step}$ depends on hyperparameters $\gamma_{u}$ and $\gamma_{d}$. It is recommended to use $\gamma_{d}=0.9$ and $\gamma_{u}=2$. This method also needs an initial choice for the Lipschitz parameter $\mu$ between $\gamma_{d}\scrL_{k-1}$ and $\scrL_{k-1}$. A choice that is reported to work well is 
\begin{align*}
\mu=\textnormal{Clip}_{[\gamma_{d}\scrL_{k-1},\scrL_{k-1}]}\left(\frac{\gap(x^{k})^{2}}{2(f(x^{k})-f(x^{k-1}))\norm{s^{k}-x^{k}}^{2}_{2}}\right).
\end{align*}

\subsection{Overhead of the backtracking}
Evaluation of the sufficient decrease condition in Algorithm 3 requires two extra evaluations of the objective function. If the condition is verified, then it is only evaluated at the current and next iterate. Following \citet{Nes13} we have the following estimate on the number of necessary function evaluations during a single execution of the backtracking procedure.

\begin{proposition}\label{prop:iterate}
Let $N_{k}$ be the number of function evaluations of the sufficient decrease condition up to iteration $k$. Then 
\begin{align*}
N_{k}\leq &(k+1)\left(1-\frac{\ln(\gamma_{d})}{\ln(\gamma_{u})}\right)\\
&+\frac{1}{\ln(\gamma_{u})}\max\{0,\ln\left(\frac{\gamma_{u} L_{\nabla f}}{\scrL_{-1}}\right)\}
\end{align*}
\end{proposition}
\begin{proof}
Call $m_{k}\geq 0$ the number of function evaluations needed in executing Algorithm 3 at stage $k$. Since the algorithm multiples the current Lipschitz parameter $\scrL_{k-1}$ by $\gamma_{u}>1$ every time that the sufficient decrease condition is not satisfied, we know that $\scrL_{k}\geq\gamma_{d}\scrL_{k-1}\gamma_{u}^{m_{k}-1}$. Hence, 
\begin{align*}
m_{k}\leq 1+\ln\left(\frac{\scrL_{k}}{\scrL_{k-1}}\right)\frac{1}{\ln(\gamma_{u})}-\frac{\ln(\gamma_{d})}{\ln(\gamma_{u})}. 
\end{align*}
Since $N_{k}=\sum_{i=0}^{k}m_{i}$, we conclude 
\begin{align*}
N_{k}\leq (k+1)\left(1-\frac{\ln(\gamma_{d})}{\ln(\gamma_{u})}\right)+\frac{1}{\ln(\gamma_{u})}\ln\left(\frac{\scrL_{k}}{\scrL_{-1}}\right).
\end{align*}
By definition of the Lipschitz parameters, we see that $\scrL_{k}\leq\max\{\gamma_{u}L_{\nabla f},\scrL_{-1}\}$. Hence, we can bound $\ln\left(\frac{\scrL_{k}}{\scrL_{-1}}\right)\leq\max\{0,\ln\left(\frac{\gamma_{u}L_{\nabla f}}{\scrL_{-1}}\right)\}$. 
\end{proof}

Proposition \ref{prop:iterate} implies that most of the backtracking subroutines terminate already after a single evaluation of the objective function gradient. Indeed, if we choose hyperparameters as $\gamma_{d}=0.9$ and $\gamma_{u}=2$, then $1-\frac{\ln(\gamma_{d})}{\ln(\gamma_{u})}\leq 1.16$ and so, asymptotically, no more than 16\% of the iterates will result in more than one gradient evaluation. 

\subsection{Proof of Theorem 3.3}
The proof of Theorem 3.4 needs the next auxiliary result which we establish first. 

\begin{lemma}\label{lem:AFWdescent}
We have for all $t\in[0,1]$
\begin{align*}
f(x^{k+1})\leq f(x^{k})-t\gap(x^{k})+\frac{t^{2}\scrL_{k}}{2}\norm{s^{k}-x^{k}}^{2}.
\end{align*}
\end{lemma}
\begin{proof}
Consider the following quadratic optimization problem 
\begin{align*}
\min_{t\in[0,1]}\{-t\gap(x^{k})+\frac{\scrL_{k}t^{2}}{2}\norm{s^{k}-x^{k}}^{2}\}.
\end{align*}
This has the unique solution 
\begin{align*}
\alpha_{k}=\tau_{k}(\scrL_{k})=\min\left\{1,\frac{\gap(x^{k})}{\scrL_{k}\norm{s^{k}-x^{k}}^{2}}\right\}.
\end{align*}
It therefore follows, 
\begin{align*}
-\alpha_{k}\gap(x^{k})&+\frac{\alpha^{2}_{k}\scrL_{k}}{2}\norm{s^{k}-x^{k}}^{2}\\
&\leq -t\gap(x^{k})+\frac{t^{2}\scrL_{k}}{2}\norm{s^{k}-x^{k}}^{2}.
\end{align*}
By definition of the backtracking procedure, Algorithm 3, we conclude
\begin{align*}
f(x^{k+1})&=f(x^{k}+\alpha_{k}(s^{k}-x^{k}))\leq Q(x^{k},\alpha_{k},\scrL_{k})\\
&=f(x^{k})-\alpha_{k}\gap(x^{k})+\frac{\alpha^{2}_{k}\scrL_{k}}{2}\norm{s^{k}-x^{k}}^{2}\\
&\leq  f(x^{k})-t\gap(x^{k})+\frac{t^{2}\scrL_{k}}{2}\norm{s^{k}-x^{k}}^{2}
\end{align*}
for all $t\in[0,1]$.
\end{proof}

\begin{proof}[Proof of Theorem 3.3.]
Define the Fenchel conjugate 
\begin{equation}
f^{\ast}(u)\eqdef\sup_{z\in\dom (f)}\{\inner{z,u}-f(z)\}.
\end{equation}
Since $f$ is proper, closed and convex, so is the Fenchel conjugate $f^{\ast}$. Moreover, since $f$ is smooth and convex on $\dom f$, we know that $f^{\ast}(u)=\inner{\nabla f(z^{\ast}(u)),u}-f(z^{\ast}(u))$, where $z^{\ast}(u)$ is the unique solution to the equation $\nabla f(z^{\ast}(u))=u$. By definition, we have $f^{\ast}(\nabla f(x^{k}))\geq \inner{\nabla f(x^{k}),x^{k}}-f(x^{k})$, and by convexity, we know that $\inner{\nabla f(x^{k}),u}-f(u)\leq \inner{\nabla f(x^{k}),x^{k}}-f(x^{k})$ for all $u\in\dom f$. We conclude, that 
\begin{align*}
f^{\ast}(\nabla f(x^{k}))\leq \inner{\nabla f(x^{k}),x^{k}}-f(x^{k}).
\end{align*}
Hence, actually equality must hold between both sides, i.e.
\begin{equation}
f^{\ast}(\nabla f(x^{k}))=\inner{\nabla f(x^{k}),x^{k}}-f(x^{k}).
\end{equation} 
Define the support function $H_{\scrX}(c)\eqdef\sup_{x\in\scrX}\inner{c,x}$, and
\begin{equation}
\psi(z)\eqdef -f^{\ast}(z)-H_{\scrX}(-z)\quad\forall z\in\dom f^{\ast}.
\end{equation}
We obtain the following series of equivalences:
\begin{align*}
\gap(x^{k})&=\inner{\nabla f(x^{k}),x^{k}-s^{k}}\\
&=\inner{\nabla f(x^{k}),x^{k}}+\inner{-\nabla f(x^{k}),s^{k}}\\
&=\inner{\nabla f(x^{k}),x^{k}}+H_{\scrX}(-\nabla f(x^{k}))\\
&=f^{\ast}(\nabla f(x^{k}))+f(x^{k})+H_{\scrX}(-\nabla f(x^{k}))\\
&=f(x^{k})-\psi(\nabla f(x^{k})).
\end{align*}
We note that $\psi$ is concave, and a dual objective function to $f$. Indeed, by the Fenchel-Young inequality, we know that $f(x)+f^{\ast}(y)\geq \inner{y,x}$ for all $x\in\dom f$ and $y\in\dom f^{\ast}$. From this inequality, we readily deduce that 
\begin{align*}
\min_{x\in\scrX}f(x)\geq -f^{\ast}(y)-H_{\scrX}(-y).
\end{align*}
Therefore, 
\begin{equation}
f^{\ast}\eqdef \min_{x\in\scrX}f(x)=\max_{y\in\dom f^{\ast}}\psi(y)\eqdef \psi^{\ast}.
\end{equation} 
We know that for all $t\in[0,1]$, 
\begin{align*}
f(x^{k+1})\leq f(x^{k})-t\gap(x^{k})+\frac{t^{2}\scrL_{k}}{2}\norm{s^{k}-x^{k}}^{2}_{2}
\end{align*}
Let us introduce the auxiliary sequence $y^{0}=\nabla f(x^{0})$, and $y^{k+1}=(1-\xi_{k})y^{k}+\xi_{k}\nabla f(x^{k})$ where $\xi_{k}\eqdef \frac{2}{k+3}$. We observe that 
\begin{align*}
f(x^{k})-\psi(y^{k})&=f(x^{k})-f^{\ast}+\psi^{\ast}-\psi(y^{k})\\
&\geq f(x^{k})-f^{\ast}.
\end{align*}
Hence, defining the approximation error $h_{k}\eqdef f(x^{k})-f^{\ast}$, we see $f(x^{k})-\psi(y^{k})\geq h_{k}$ for all $k\geq 0$. Moreover, since $\psi$ is concave, we know that 
\begin{align*}
\psi(y^{k+1})\geq (1-\xi_{k})\psi(y^{k})+\xi_{k}\psi(\nabla f(x^{k})).
\end{align*}
Consequently, 
\begin{align*}
h_{k+1}&\leq f(x^{k+1})-\psi(y^{k+1})\\
&\leq f(x^{k})-\xi_{k}\gap(x^{k})+\frac{\xi^{2}_{k}\scrL_{k}}{2}\norm{s^{k}-x^{k}}^{2}_{2}\\
&-(1-\xi_{k})\psi(y^{k})-\xi_{k}\psi(\nabla f(x^{k}))\\
&=(1-\xi_{k})[f(x^{k})-\psi(y^{k})]+\frac{\xi^{2}_{k}\scrL_{k}}{2}\norm{s^{k}-x^{k}}^{2}_{2}\\
&\leq (1-\xi_{k})[f(x^{k})-\psi(y^{k})]+\frac{\xi^{2}_{k}\scrL_{k}}{2}\diam(\scrX)^{2}.
\end{align*}
Define $A_{k}\eqdef \frac{1}{2}(k+1)(k+2)$ for $k\geq 0$. For this specification, it is easy to check that 
\begin{align}
&A_{k+1}(1-\xi_{k})=A_{k},\text{ and }\\
&A_{k+1}\frac{\xi^{2}_{k}}{2}\leq 1.
\end{align}
Hence, 
\begin{align*}
A_{k+1}[f(x^{k+1})-\psi(y^{k+1})]&\leq A_{k+1}(1-\xi_{k})[f(x^{k})-\psi(y^{k})]\\
&+A_{k+1}\frac{\xi^{2}_{k}}{2}\scrL_{k}\diam(\scrX)^{2}\\
&\leq A_{k}[f(x^{k})-\psi(y^{k})]\\
&+\scrL_{k}\diam(\scrX)^{2}.
\end{align*}
Summing from $i=0,\ldots,k-1$, and calling 
$$
\bar{\scrL}_{k}\eqdef\frac{1}{k}\sum_{i=0}^{k-1}\scrL_{i},
$$
 this implies 
\begin{align*}
&f(x^{k})-\psi(y^{k})\leq\frac{1}{A_{k}}[f(x^{0})-\psi(y^{0})]+\frac{k\diam(\scrX)^{2}}{2A_{k}}\bar{\scrL}_{k}\\
&=\frac{2}{(k+1)(k+2)}[f(x^{0})-\psi(y^{0})]+\frac{k\diam(\scrX)^{2}}{(k+1)(k+2)}\bar{\scrL}_{k}\\
&=\frac{2\gap(x^{0})}{(k+1)(k+2)}+\frac{k\diam(\scrX)^{2}}{(k+1)(k+2)}\bar{\scrL}_{k}\\
\end{align*}
Since $\scrL_{k}\leq L_{\nabla f}$, we get $h_{k}=\scrO(k^{-1})$. 
\end{proof}

\section{Proof of Theorem 4.2}
\label{sec:LLOO}
Let us define $\scrP(x^{0})\eqdef \left\{x\in\scrX: f(x)\leq f^{\ast}+\gap(x^{0})\right\}$. We prove this theorem by induction. For $k=0$, we have the given initial condition $x^{0}\in\dom f\cap\scrX$. Since $x^{0}\in\scrP(x^{0})$ trivially, we know from Lemma \ref{lem:gbound} that 
\begin{align*}
h_{0}\geq\frac{\sigma_{f}}{6}\norm{x^0-x^{\ast}}^{2}_{2}.
\end{align*}
Set $r_{0}\geq\sqrt{\frac{6h_{0}}{\sigma_{f}}}$, this implies $x^{\ast}\in B(x^{0},r_{0})$. Since $s^{0}=\scrA(x^{0},r_{0},\nabla f(x^{0}))$, the definition of the \ac{LLOO} tells us 
\begin{align*}
f(x^{\ast})-f(x^{0})&\geq \inner{\nabla f(x^{0}),x^{\ast}-x^{0}}\\
&\geq\inner{\nabla f(x^{0}),s^{0}-x^{0}}\geq-\gap(x^0),
\end{align*}
where the first inequality is a consequence of the convexity of $f$, whereas the last inequality follows by definition of the dual gap function. Therefore, we have $h_{0}\leq \gap(x^{0})$. Set $r_{0}=\sqrt{\frac{6\gap(x^{0})}{\sigma_{f}}}$, and let $\alpha_0=\min\left\{\frac{\gap(x^{0})}{\frac{4}{M^{2}}(\ce^0)^2},1\right\}\frac{1}{\ce^0+1}<1$.
Note that this choice of $\alpha_0$ guarantees that $\alpha_0\ce^0\leq \frac{\ce^0}{\ce^0+1}<1$. Hence, doing the update $x^{1}=x^{0}+\alpha_{0}(s^{0}-x^{0})$ we know from convexity of $\scrX$ and Lemma \ref{lem:Dikin} that $x^{1}\in\dom f\cap\scrX$. 

Apply inequality \eqref{eq:up} to conclude
\begin{align*}
f(x^{1})\leq f(x^{0})+\alpha_{0}\inner{\nabla f(x^{0}),s^{0}-x^{0}}+\frac{4}{M^{2}}\omega_{\ast}(\alpha_{0}\ce^{0}).
\end{align*}
Since, $\omega_{\ast}(t)=-t-\ln(1-t)$, it follows 
 \begin{align*}
 \omega_{\ast}(t)=\sum_{j=2}^{\infty}\frac{t^{j}}{j}\leq\frac{t^{2}}{2}\sum_{j=0}^{\infty}t^{j}=\frac{t^{2}}{2(1-t)}\quad\forall t\in[0,1).
 \end{align*}
 Since $\alpha_{0}\ce^{0}\in[0,1)$, for all $k\geq 0$, we therefore arrive at the estimate
 \begin{align*}
h_{1}\leq h_{0}+\alpha_{0}(f(x^{\ast})-f(x^{0}))+\frac{4}{M^{2}}\frac{\alpha_{0}^{2}(\ce^{0})^{2}}{2(1-\alpha_{0}\ce^{0})}.
\end{align*}
By the above said, we know that $1-\alpha_{0}\ce^{0}\geq \frac{1}{1+\ce^{0}}$, and therefore, 
\begin{align*}
h_{1}&\leq (1-\alpha_{0})h_{0}+\frac{2\alpha_{0}^{2}}{M^{2}}(\ce^{0})^{2}(1+\ce^0)\\
&\leq (1-\alpha_{0})h_{0}+\frac{2\alpha_{0}^{2}}{M^{2}}(\ce^{0})^{2}(1+\ce^0)\\
&\leq  (1-\alpha_{0})\gap(x^{0})+\frac{2\alpha_{0}^{2}}{M^{2}}(\ce^{0})^{2}(1+\ce^0)
\end{align*}
Plugging in the chosen value of $\alpha_0$ we obtain that
\begin{align*}h_1&\leq \gap(x^0)+\alpha_0(-\gap(x^0)+\frac{2\alpha_{0}}{M^{2}}(\ce^{0})^{2}(1+\ce^0))\\
&\leq \gap(x^0)(1-\frac{\alpha_0}{2}) 
\end{align*}
Notice that by the definition of the LLOO we have that
\begin{align*}
\ce^0&\leq \sqrt{L_{\nabla f}}\frac{M}{2}\norm{x^0-s^0}\\
&\leq \frac{\sqrt{L_{\nabla f}}M}{2}\min\{\rho r_{0},\diam(\scrX)\}
\end{align*}
which implies that 
\begin{align*}\frac{\gap(x^0)}{\frac{4}{M^{2}}(\ce^0)^2}\geq \frac{\gap(x^0)}{L_{\nabla f} \rho^2 r_0^2}\geq \frac{\sigma_f}{6L_{\nabla f} \rho^2}\end{align*}
where the last inequality follows from the definition of $r_0$.

Thus, we have that
\begin{align*}
\alpha_0\geq \min\{ \frac{\sigma_f}{6L_{\nabla f} \rho^2},1\}\frac{1}{1+ \sqrt{L_{\nabla f}}M\diam(\scrX)/2}\equiv \bar{\alpha},
\end{align*}
which implies that 
\begin{align*}
h_1&\leq \gap(x^0)(1-\frac{\alpha_0}{2})\\
&\leq \gap(x^0)\exp(-\frac{\alpha_{0}}{2})\\
&=\gap(x^{0})c_{1}.
\end{align*}
This verifies the claim for $k=0$. Now proceed inductively. Suppose that 
\begin{equation}\label{eq:IH}
h_{k}\leq\gap(x^{0})c_{k},\, c_{k}\eqdef\exp\left(-\frac{1}{2}\sum_{i=0}^{k-1}\alpha_{i}\right).
\end{equation}
Then $x^{k}\in\scrP(x^{0})$, and Lemma \ref{lem:gbound} tells us that 
\begin{align*}
\norm{x^{k}-x^{\ast}}^{2}_{2}\leq \frac{6 h_{k}}{\sigma_{f}}\leq \frac{6\gap(x^0)}{\sigma_f}c_k=r_0^2c_k\equiv r^{2}_{k}.
\end{align*}
Hence, $x^{\ast}\in\ball(x^{k},r_{k})$. Proceeding as for $k=0$, let us again make the educated guess that we can take a step size $\alpha_{k}<1/c_{k}$. By the same argument as before, we obtain sufficient decrease 
\begin{align*}
h_{k+1}\leq h_{k}+\alpha_{k}\inner{\nabla f(x^{k}),s^{k}-x^{k}}+\frac{4}{M^{2}}\frac{\alpha^{2}_{k}(\ce^{k})^{2}}{2(1-\alpha_{k}\ce^{k})}.
\end{align*}
Since $s^{k}=\scrA(x^{k},r_{k},\nabla f(x^{k}))$, we know by definition of the \ac{LLOO} that $\inner{\nabla f(x^{k}),x^{\ast}-x^{k}}\geq \inner{\nabla f(x^{k}),s^{k}-x^{k}}$ and $\norm{s^{k}-x^{k}}_{2}\leq\rho r_{k}$. Consequently,  setting $\alpha_k=\min\{\frac{c_k\gap(x^0)}{\frac{4}{M^{2}}(\ce^k)^2},1\}\frac{1}{(1+\ce^k)}$ we obtain that
\begin{align*}
h_{k+1}&\leq(1-\alpha_{k})h_{k}+\frac{2}{M^{2}}\frac{\alpha^{2}_{k}(\ce^{k})^{2}}{1-\alpha_{k}\ce^{k}}\\
&\leq (1-\alpha_{k})\gap(x^0)c_k+\frac{2}{M^{2}}\alpha^{2}_{k}(\ce^{k})^{2}(1+\ce^{k})\\
&\leq \gap(x^0)c_k(1-\frac{\alpha_{k}}{2})
\end{align*}
where the second inequality follow from the fact that
\begin{align*}
&\alpha_k\ce^k
\leq\frac{\ce^k}{1+\ce^k}<1,\text{ and }\\
&1-\alpha_{k}\ce^{k}\geq\frac{1}{1+\ce^{k}}. 
\end{align*}
as well as the upper bound on $h_k$ obtained by the induction step.
Finally by the definition of the \ac{LLOO} we have that
\begin{align*}
\ce^k&\leq \sqrt{L_{\nabla f}}\frac{M}{2}\norm{x^k-s^k}\\
&\leq \frac{\sqrt{L_{\nabla f}}M}{2}\min\{\rho r_k, \diam(\scrX)\},
\end{align*}
which implies that 
\begin{align*}
\frac{M^2\gap(x^0)c_k}{4(\ce^k)^2}\geq \frac{\gap(x^0)c_k}{L_{\nabla f} \rho^2 r_k^2}\geq \frac{\sigma_f}{6L \rho^2}
\end{align*}
where the last inequality follows from the definition of $r_k$.
Thus, we have that
\begin{align*}
\alpha_k\geq \min\{ \frac{\sigma_f}{6L_{\nabla f} \rho^2},1\}\frac{1}{1+ \sqrt{L_{\nabla f}}\frac{M\diam(\scrX)}{2}}\equiv \bar{\alpha},
\end{align*}
which implies that 
\begin{align*}
h_{k+1}&\leq\gap(x^0)c_k(1-\frac{\alpha_k}{2})\\
&\leq \gap(x^0)c_k\exp\left(-\frac{\alpha_{k}}{2}\right)\\
&=\gap(x^{0})c_{k+1}.
\end{align*}
Since $\alpha_{k}\geq\bar{\alpha}$ for all $k$, we thus have shown that 
\begin{equation}
h_{k}\leq \gap(x^{0})\exp(-k\bar{\alpha}/2).
\end{equation}

\section{Numerical Experiments}
\label{sec:numerics}
We give extensive information about the numerical experiments we have conducted in order to test the performance of all three algorithms developed in this paper.  In the numerical experiments we tested the performance of Variant 1 (V1) and Variant 2 (V2) of Algorithm 2, and compared them with the performance of Frank-Wolfe with {standard} step-size of $\frac{2}{k+2}$ (Standard), and step-size determined by exact {line-search} (Line-S.). As a further benchmark, the self-concordant Proximal-Newton (PN) of \citet{CevKyrTra15}, as implemented in the SCOPT package\footnote{\url{https://www.epfl.ch/labs/lions/technology/scopt/}}, is included. For the portfolio optimization problem, Algorithm 4 is also implemented. All codes are written in Python 3, with packages for scientific computing NumPy 1.18.1 and SciPy 1.4.1. The experiments were conducted on a PC with Intel Core i5-7500 3.4GHzs, with a total of 16GB RAM.

In both experiments the Frank-Wolfe based methods have been terminated after 50,000 iterations. Because of its higher computational complexity, we decided to stop PN after 1,000 iterations. Each algorithm was terminated early if the optimality gap in a given iteration was lower than $1e-10$. LLOO was only implemented for the portfolio selection problem, using the local linear oracle given in \cite{GarHaz16} for the simplex, as described in Appendix F in the supplementary materials.

\subsection{Results on the Portfolio Optimization problem}
For the Portfolio Optimization problem we used synthetic data, as in Section 6.4 in \cite{SunTran18}. The details of the data generating process are as follows. We generate matrix $R$ with given price ratios as: $R := ones(n, p) + N (0, 0.1)$, which allows the closing price to vary about 10\% between two consecutive periods. We used different sizes of matrix $R$: $(n, p) = (1000, 800), (1000, 1200)$, and $(1000, 1500)$ with 4 samples for each size. Hence, there are totally 12 datasets. The detailed results for 9 out of these 12 datasets is reported in Figure \ref{fig:Portfolio}.

\begin{figure*}
\label{fig:Portfolio}
\centering
\begin{tabular}{|c|c|c|}
\hline
\subf{\includegraphics[width=50mm]{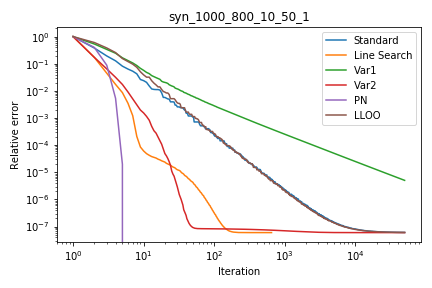}}
     {}
&
\subf{\includegraphics[width=50mm]{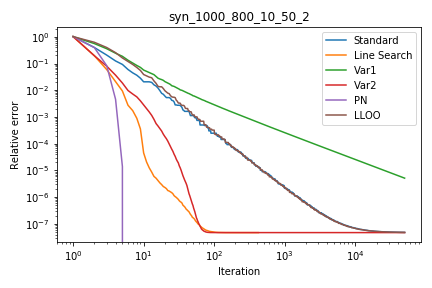}}
     {}
&
\subf{\includegraphics[width=50mm]{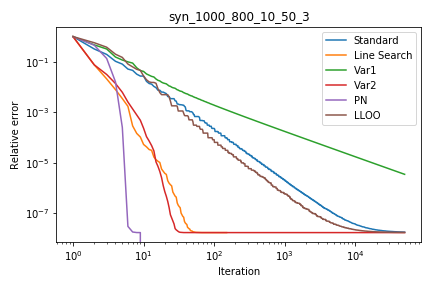}}
     {}
\\
\hline
\subf{\includegraphics[width=50mm]{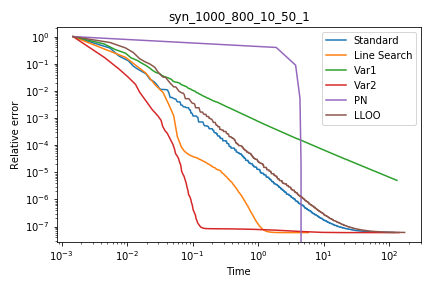}}
     {}
&
\subf{\includegraphics[width=50mm]{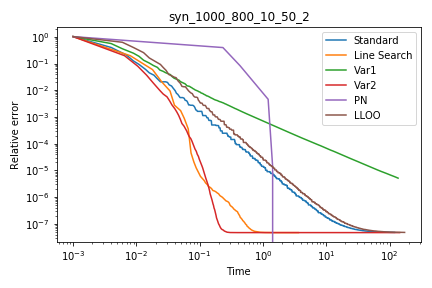}}
     {}
&
\subf{\includegraphics[width=50mm]{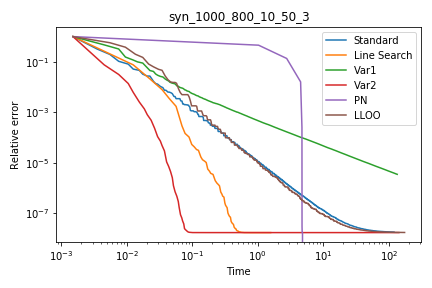}}
     {}
\\     
\hline
\subf{\includegraphics[width=50mm]{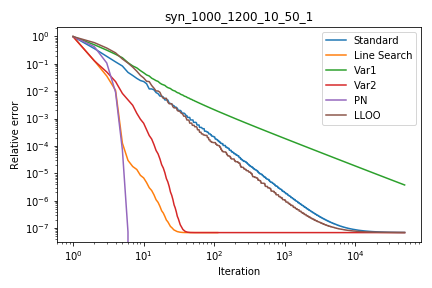}}
     {}
&
\subf{\includegraphics[width=50mm]{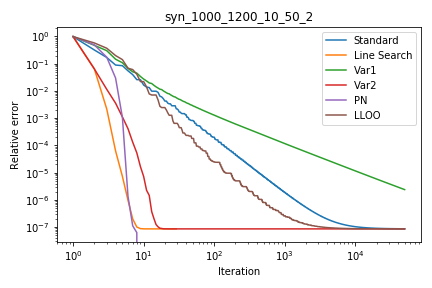}}
     {}
&
\subf{\includegraphics[width=50mm]{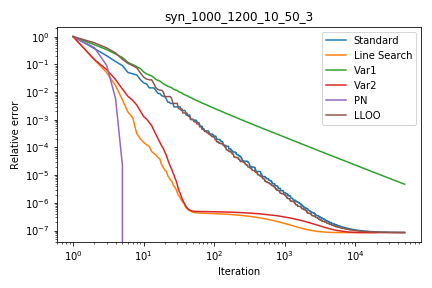}}
     {}     
\\
\hline
\subf{\includegraphics[width=50mm]{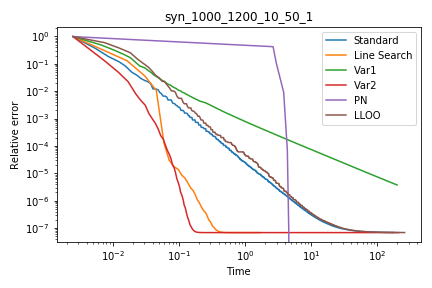}}
     {}
&
\subf{\includegraphics[width=50mm]{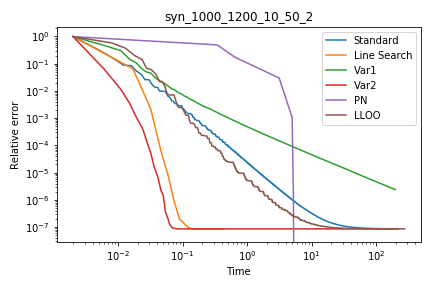}}
     {}
&
\subf{\includegraphics[width=50mm]{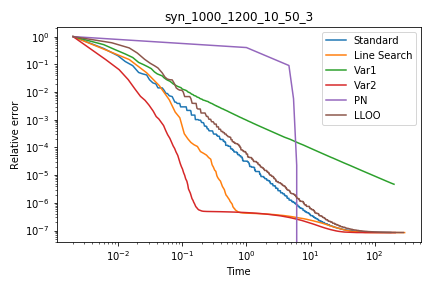}}
     {}
\\
\hline
\subf{\includegraphics[width=50mm]{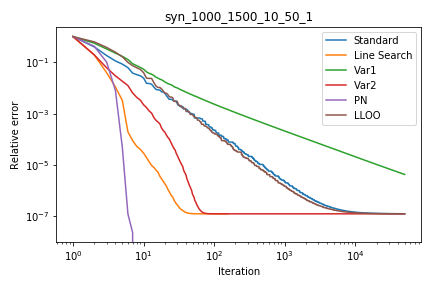}}
     {}
&
\subf{\includegraphics[width=50mm]{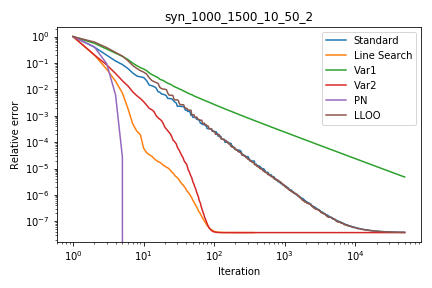}}
     {}
&
\subf{\includegraphics[width=50mm]{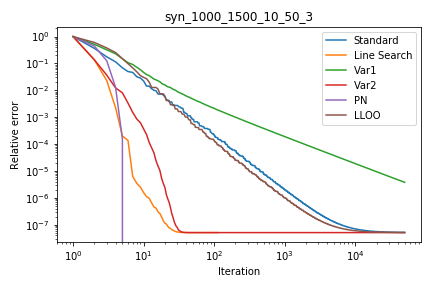}}
     {}
\\
\subf{\includegraphics[width=50mm]{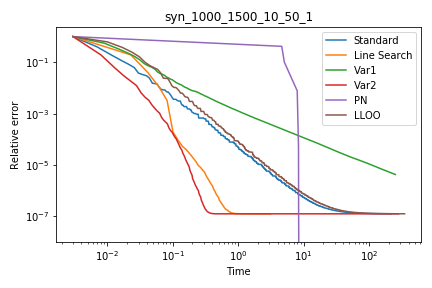}}
     {}
&
     \subf{\includegraphics[width=50mm]{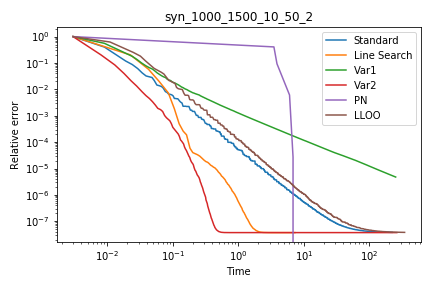}}
        {}
   &
     \subf{\includegraphics[width=50mm]{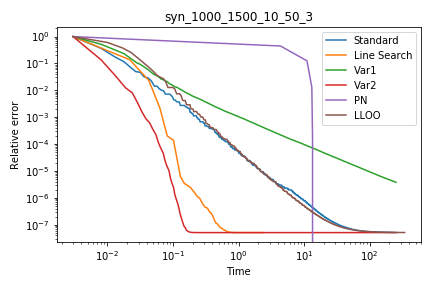}}
     {}
\\
\hline
\end{tabular}
\end{figure*}

\subsection{Results on the Poisson Inverse problem}
\begin{figure*}
\label{fig:poisson}
\centering
\begin{tabular}{|c|c|c|}
\hline
\subf{\includegraphics[width=50mm]{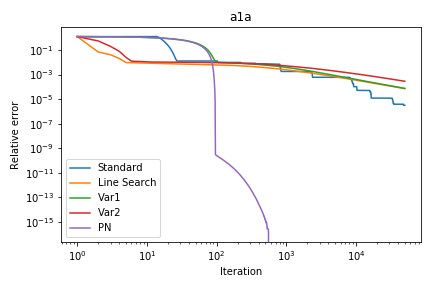}}
    {}
&
\subf{\includegraphics[width=50mm]{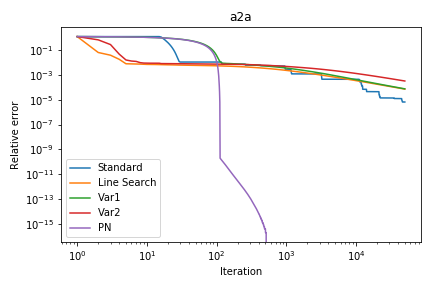}}
     {}
&
\subf{\includegraphics[width=50mm]{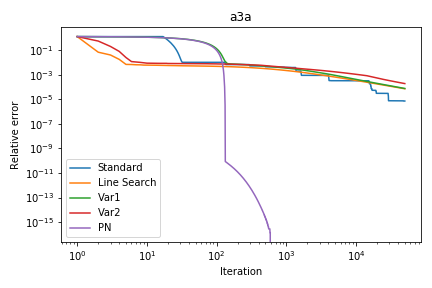}}
     {}
\\
\hline
\subf{\includegraphics[width=50mm]{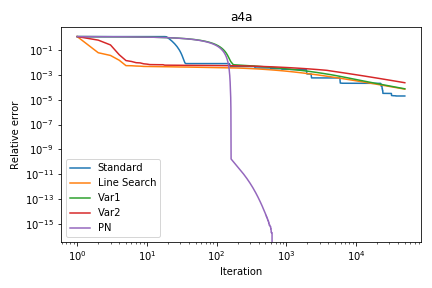}}
     {}
&
\subf{\includegraphics[width=50mm]{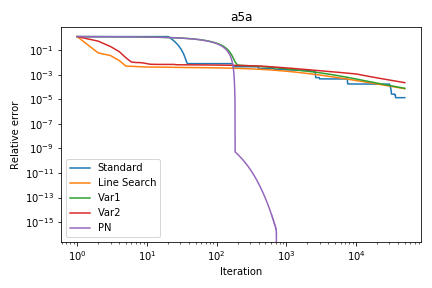}}
    {}
&
\subf{\includegraphics[width=50mm]{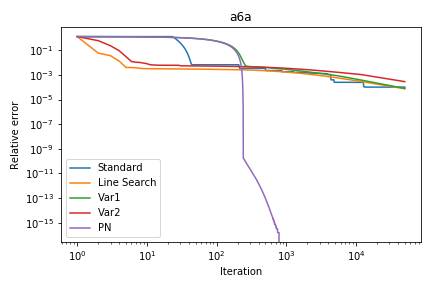}}
     {}
\\     
\hline
\subf{\includegraphics[width=50mm]{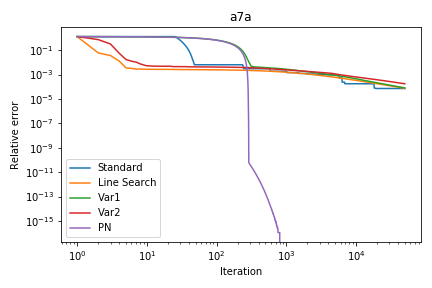}}
     {}
&
\subf{\includegraphics[width=50mm]{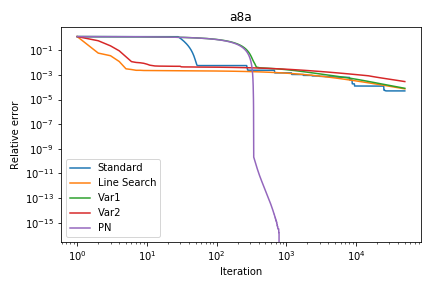}}
     {}
&
\subf{\includegraphics[width=50mm]{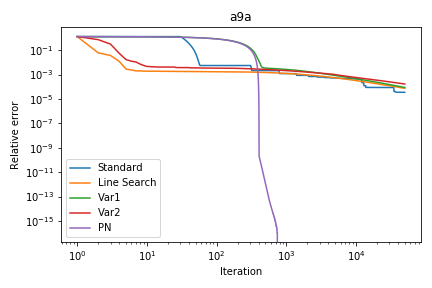}}
     {}     
\\
\hline
\subf{\includegraphics[width=50mm]{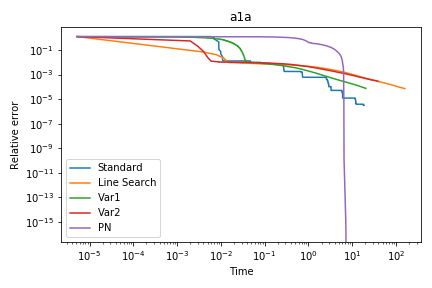}}
     {}
&
\subf{\includegraphics[width=50mm]{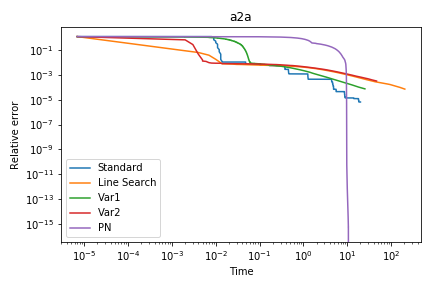}}
     {}
&
\subf{\includegraphics[width=50mm]{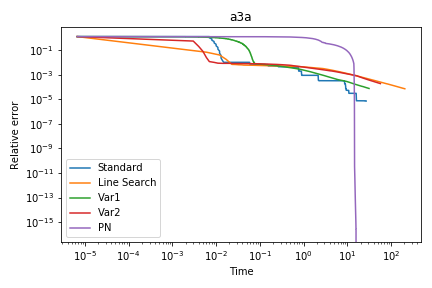}}
     {}
\\
\hline
\subf{\includegraphics[width=50mm]{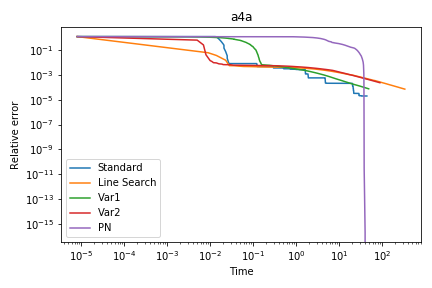}}
     {}
&
\subf{\includegraphics[width=50mm]{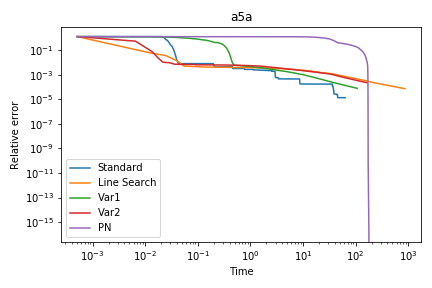}}
     {}
&
\subf{\includegraphics[width=50mm]{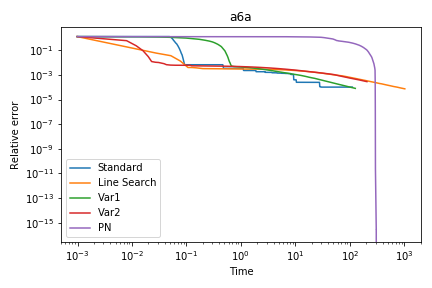}}
     {}
\\
\subf{\includegraphics[width=50mm]{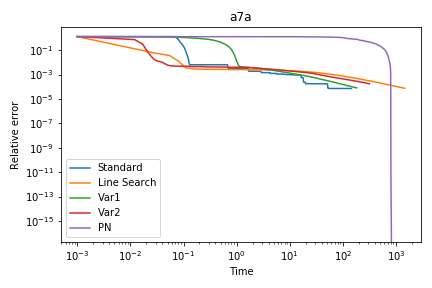}}
     {}
&
\subf{\includegraphics[width=50mm]{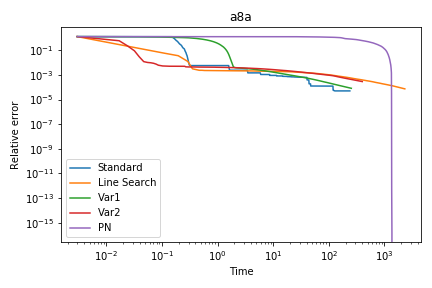}}
        {}
   &
   \subf{\includegraphics[width=50mm]{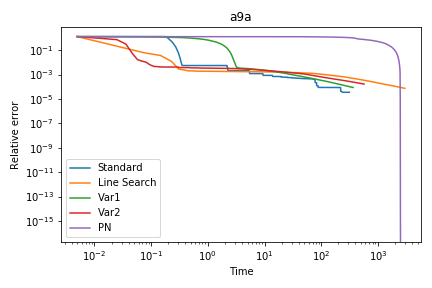}}
   {}
\\
\hline
\end{tabular}
\end{figure*}

For the Poisson inverse problem we used the datasets a1a-a9a from the LIBSVM library \cite{LIBSVM}. The results on each individual data set are displayed in Figure \ref{fig:poisson}. 

\section{Constructing an LLOO for Simplex Constraints}
\label{sec:Simplex}
In this section we outline a construction of an LLOO in case where the polytope is the unit simplex $\scrX\equiv\Delta_{n}=\{x\in\Rn\vert x\geq 0,\sum_{i=1}^{n}x_{i}=1\}$. This construction is given in \cite{GarHaz16}. 
\begin{lemma}
Given a point $x\in\Delta_{n}$, a radius $r>0$, and a linear objective $c\in\Rn$, consider the optimization problem
\begin{equation}\label{eq:l1}
\min\{\inner{c,y}\vert\; \norm{x-y}_{1}\leq d\}
\end{equation}
for some $d > 0$. Let us denote by $p^{\ast}$ an optimal solution to problem \eqref{eq:l1} when we set $d=\sqrt{n}r$. Then $p^{\ast}$ is the output of an LLOO with parameter $\rho=\sqrt{n}$ for $\scrX=\Delta_{n}$. That is, for all $y\in\Delta_{n}\cap B(x,r)$
\begin{equation}
\inner{p^{\ast},c}\leq \inner{y,c}\text{ and }\norm{x-p^{\ast}}\leq\sqrt{n}r.
\end{equation}
\end{lemma}

Algorithm \ref{alg:LLOOsimplex} implements an LLOO for the unit simplex. In this algorithm we use the Kronecker delta $\delta_{i,j}=1$ if $j=i$ and $0$ otherwise. 
\begin{algorithm}
 \caption{LLOO for the simplex}
 \label{alg:LLOOsimplex}
 \begin{algorithmic}
\STATE {\bfseries Input: }
$x\in\Delta_{n}$, radius $r>0$, cost vector $c\in\Rn$.
\STATE Set $d=\sqrt{n}r$
\STATE $m=\min\{d/2,1\}$
\STATE $i^{\ast}=\argmin_{1\leq i\leq n}c_{i}$
\STATE $p_{+}=m[\delta_{i^{\ast},1};\ldots;\delta_{i^{\ast},n}]^{\top}$.
\STATE $p_{-}=\0_{n}\in\Rn$
\STATE Let $i_{1},\ldots,i_{n}$ be a permutation over $\{1,\ldots,n\}$ such that $c_{i_{1}}\geq c_{i_{2}}\geq\ldots\geq c_{i_{n}}$.
\STATE Set $k=\min\{\ell \vert \sum_{j=1}^{\ell}x_{i_{j}}\geq m\}$
\STATE For all $1\leq j\leq k-1$ set $(p_{-})_{i_{j}}=x_{i_{j}}$
\STATE $p_{i_{k}}=m-\sum_{j=1}^{k-1}x_{i_{j}}$ 
\STATE Update $p=x+p_{+}-p_{-}$.
\end{algorithmic}
\end{algorithm}

\end{document}